\documentclass[12pt]{article}

\usepackage{amsfonts}
\usepackage{latexsym}
\usepackage{pgf,tikz}
\usetikzlibrary{arrows}
\usetikzlibrary[patterns]

\usepackage{setspace}

\usepackage{graphicx} % Required for including images
\graphicspath{{Figures/}} % Set the default folder for images

\usepackage{enumitem} % Required for manipulating the whitespace between and within lists

\usepackage{lipsum} % Used for inserting dummy 'Lorem ipsum' text into the template

\usepackage{subfig} % Required for creating figures with multiple parts (subfigures)

\usepackage{amsmath,amssymb,amsthm} % For including math equations, theorems, symbols, etc

\usepackage{varioref} % More descriptive referencing

\usepackage{bm}
\usepackage{bbold}
\usepackage[makeroom]{cancel}
\usepackage{lscape}
\usepackage{mathtools}
\usepackage{caption}
\usepackage{multicol}
\usepackage{tikz}
\usepackage{float}

\usepackage[toc,title]{appendix}

\usepackage[text={6.5in,9in}]{geometry}
\theoremstyle{remark} % Define theorem styles here based on the remark style (used for remarks and notes)

\newcommand {\be}{\begin{equation}} 
\newcommand {\ee}{\end{equation}}
\newcommand{\tr}{\mathrm{tr}}
\newcommand{\intO}[1]{\int_{\Omega}#1 \;dx}
\newcommand{\intG}[2]{\int_{\Gamma_#1}#2 \;ds}
\newcommand{\divg}{\mathrm{div}\;}
\newcommand{\Vnorm}[1]{||#1||_{\mathcal{V}}}
\newcommand{\Znorm}[1]{||#1||_{0}}

\begin{document}

\begin{center}
{\Large {\bf Finite element approximations for near-incompressible and near-inextensible transversely isotropic bodies}}
\vspace{3ex}\\
{\large F. Rasolofoson \quad B.J. Grieshaber \quad B.D. Reddy}
\vspace{2ex}\\
Department of Mathematics and Applied Mathematics, \\
and Centre for Research in Computational and Applied Mechanics, University of Cape Town, South Africa
\vspace{1ex}\\
(rslfar002@myuct.ac.za,beverley.grieshaber@uct.ac.za,daya.reddy@uct.ac.za)
\end{center}

\section*{Abstract}
This work comprises a detailed theoretical and computational study of the boundary value problem for transversely isotropic linear elastic bodies.
General conditions for well-posedness are derived in terms of the material parameters.
The discrete form of the displacement problem is formulated for conforming finite element approximations.
The error estimate reveals that anisotropy can play a role in minimizing or even eliminating locking behaviour, for moderate values of the ratio of Young's moduli in the fibre and transverse directions.
In addition to the standard conforming approximation an alternative formulation, involving under-integration of the volumetric and extensional terms in the weak formulation, is considered.
The latter is equivalent to either a mixed or a perturbed Lagrangian formulation, analogously to the well-known situation for the volumetric term.
A set of numerical examples confirms the locking-free behaviour in the near-incompressible limit of the standard formulation with moderate anisotropy, with locking behaviour being clearly evident in the case of near-inextensibility.
On the other hand, under-integration of the extensional term leads to extensional locking-free behaviour, with convergence at superlinear rates.

{\bf Keywords:} solids, finite element methods, elasticity

\section{Introduction}
Anisotropy is a significant mechanical feature of composite materials, for example, in the aerospace and automotive industries.
It occurs naturally in crystalline structures, geotechnical materials, and in  the mechanical properties of biological media such as muscles, tendons, or bones (see for example \cite{Altenbach2005,Exadaktylos2001,Humphrey2002,Ogden1984}).
Fibrous or fibre-reinforced structures are generally modelled as transversely isotropic materials, with the plane normal to the fibre direction corresponding to a plane of isotropy.
Important early contributions to the mechanics of fibre-reinforced materials include the works \cite{Hashin1970,Pipkin1979,Spencer1972}.
The monograph \cite{Ting1996} gives a detailed presentation of the mechanics of anisotropic materials.

Early work on the linear problem, in addition to some of the works cited above, includes the investigations \cite{Hayes-Horgan1974,Hayes-Horgan1975} of uniqueness and local (pointwise) stability for materials with the constraint of inextensibility.
There have been extensive theoretical and computational studies carried out of the behaviour of transversely isotropic materials, or the closely related topic of inextensible or near-inextensible materials, in the context of hyperelasticity
(see for example the key works \cite{Schroeder-Neff2003,Schroeder-Neff-Balzani2005,Weiss-Maker-Govindjee1996}).
Mixed finite element approximations leading to stable formulations have been developed and implemented for near-incompressibility and near-inextensibility in \cite{Zdunek-Rachowicz-Eriksson2016}, while transverse anisotropy has been studied, also using various mixed finite element formulations, in \cite{Zdunek-Rachowicz2017a,Zdunek-Rachowicz2017b}.
A mixed approach allied with a Lagrange multiplier or perturbed Lagrange multiplier formulation of the inextensibility constraint forms the basis of a computational study in \cite{Wriggers-Schroeder-Auricchio2014}.
In earlier work, a corresponding numerical study, using a mixed three-field approach for the linear problem, has been carried out in \cite{Bradaigh-Pipes1992}.
The recent contribution \cite{Wriggers-Hudobivnik-Korelc2017} presents a locking-free implementation for nonlinear transversely isotropic elasticity, the inextensibility constraint in the fibre direction being accommodated using a perturbed Lagrangian formulation.
A corresponding treatment of the small-strain problem has been the subject of \cite{Auricchio-Scalet-Wriggers2017}, in which limiting extensibility is investigated numerically using penalty, Lagrange multiplier, and perturbed Lagrangian approaches. 

A model of a bi-directional elastic composite has been developed in \cite{Zeidi-Kim2018}, and studied computationally for the case of plane deformations, as an extension of earlier work \cite{Zeidi-Kim2017} on single families of fibres.
The work \cite{Steigmann2012} extends conventional studies of transverse isotropy by developing a model in the context of nonlinear elasticity, in which the fibres are assumed to be resistant to flexure and twist, in addition to extension. 

There has been little work on the well-posedness of problems with internal constraints such as inextensibility, in contrast to the many treatments of incompressibility.
The investigation \cite{Arnold-Falk1987} approaches  the problem via a mixed formulation of Hellinger-Reissner type, and establishes conditions on the elastic constants for the problem to be well-posed  for inextensible materials, and also for the case of orthotropy.
Corresponding abstract results have been presented in \cite{Dantas1999} for the case of non-homogeneous materials.
In \cite{Exadaktylos2001} conditions for pointwise stability are obtained for materials in plane strain.

The purpose of this work is to undertake a detailed theoretical and computational study of the behaviour of transversely isotropic linear elastic materials, with a view to establishing conditions under which conforming finite element approximations are uniformly convergent in the incompressible and inextensible limits.
In this sense the present investigation extends considerably the work reported in \cite{Auricchio-Scalet-Wriggers2017}, by first establishing conditions for well-posedness of the continuous problem, then by establishing error estimates for conforming finite element approximations that shed light on the conditions under which locking-free behaviour can be expected.
It is shown that the presence of transverse isotropy leads to behaviour that is volumetrically locking-free, in circumstances that would lead to locking behaviour for isotropic materials.
In particular, in situations of moderate anisotropy, when the ratio of Young's modulus in the fibre direction to that in the plane of isotropy is not too large, in a sense that will be made precise, volumetric locking-free behaviour results for low-order elements.

A second contribution relates to the circumstances under which extensional locking-free behaviour occurs for low-order elements.
While the standard approximations do in fact exhibit locking, it is shown that, analogous to the now classical result for near-incompressibility \cite{Hughes1987}, under-integration of the extensional term results in locking-free behaviour.
This approach is shown to be equivalent to a mixed formulation, as in the case of incompressible behaviour, and also to a perturbed Lagrangian approach of the kind presented, for example, in \cite{Auricchio-Scalet-Wriggers2017}.

The structure of the rest of this work is as follows.
Section 2 sets out the details of constitutive relations for transversely isotropic linear elastic materials, and establishes conditions on the material constants for pointwise stability.
The weak form of the boundary-value problem is presented in Section 3, and conditions for well-posedness derived. Section 4 is concerned with conforming finite element approximations.
The standard error estimate reveals the role played by anisotropy in mitigating locking behaviour in the incompressible limit.
An alternative approximation, employing under-integration of the volumetric and extensional terms, is related to a mixed as well as a perturbed Lagrangian formulation.
These various features are explored numerically in Section 5, through two sets of example problems. The work concludes with a summary of results and a discussion of open problems.

%------------------------------------------------------------------------------------------------
%	TI materials
%------------------------------------------------------------------------------------------------

\section{Transversely isotropic materials}
For a transversely isotropic linearly elastic material with fibre direction given by the unit vector $\bm{a}$, the elasticity tensor is given by \cite{Lubarda-Chen2008,Spencer1982}
\be\label{tangent_moduli}
\mathbb{C} = \lambda \bm{I}\otimes\bm{I} + 2\mu_t \mathbb{I} + \beta \bm{M}\otimes\bm{M} + \alpha (\bm{I}\otimes\bm{M}+\bm{M}\otimes\bm{I}) + \gamma \mathbb{M}\,.
\ee
Here $\bm{I}$ is the second-order identity tensor, $\mathbb{I}$ is the fourth-order identity tensor, $\bm{M} = \bm{a}\otimes \bm{a}$, and $\mathbb{M}$ is the fourth-order tensor defined by
\be
\mathbb{M}\bm{R} =   \bm{M}\bm{R} + \bm{R}\bm{M}\qquad\mbox{for any second-order tensor}\ \bm{R}\,. 
\label{mathcalM}
\ee
$\lambda$ denotes the first Lam\'e parameter, the shear modulus in the plane of isotropy is $\mu_t$, $\mu_l$ is the shear modulus along the fibre direction, and 
\be
\gamma = 2(\mu_l - \mu_t).
\label{gamma}
\ee
The further material constants $\alpha$ and $\beta$ do not have a direct interpretation, though it will be seen that $\beta \rightarrow \infty$ in the inextensible limit.

In component form, equation \eqref{tangent_moduli} becomes
\begin{align}
\mathbb{C}_{ABCD} =& \lambda \delta_{AB}\delta_{CD}
					+ 2\mu_t\delta_{AC}\delta_{BD}
					+ \beta M_{AB}M_{CD}
					+ \alpha(\delta_{AB}M_{CD} + M_{AB}\delta_{CD})\nonumber\\
				  & + 2(\mu_l - \mu_t)(\delta_{AC}M_{BD} + M_{AC}\delta_{BD}).\label{altC_ABCD}
\end{align}
The corresponding linear stress-strain relation for small deformations is then 
\begin{align}
	\bm\sigma &= \mathbb{C}\bm{\varepsilon}\nonumber\\
			  &= \lambda(\tr\bm\varepsilon)\bm{I} + 2\mu_t \bm\varepsilon + \beta (\bm{M}:\bm\varepsilon)\bm{M} + \alpha ((\bm{M}:\bm\varepsilon)\bm{I} + (\tr\bm\varepsilon)\bm{M}) + \gamma (\bm\varepsilon\bm{M} + \bm{M}\bm\varepsilon),\label{cauchy_tensor}
\end{align}
in which $\bm\sigma$ and $\bm\varepsilon$ denote the stress and the infinitesimal strain tensors; $\mbox{tr}\,\bm{\varepsilon}$ denotes the trace of $\bm{\varepsilon}$, and $\bm{M}:\bm{\varepsilon} = \bm{\varepsilon}\bm{a}\cdot\bm{a}$, obtained from the definition of $\bm{M}$, gives the strain in the direction of $\bm{a}.$
The special case of an isotropic material is recovered by setting $\alpha = \beta = 0$ and $\mu_l = \mu_t$.

For the particular case in which $\bm{a} = \bm{e}_3$, the stress-strain relationship can be written in matrix form as 
\be\label{stress_strain_e3}
\begin{pmatrix}
	\sigma_{11} \\ \sigma_{22} \\ \sigma_{33} \\ \sigma_{23} \\ \sigma_{13} \\ \sigma_{12}
\end{pmatrix}
=
\begin{pmatrix}
	\lambda + 2\mu_t & \lambda & \lambda + \alpha & 0 & 0 & 0\\
	\lambda & \lambda + 2\mu_t & \lambda + \alpha & 0 & 0 & 0\\
	\lambda + \alpha & \lambda + \alpha & \lambda+2\mu_t+\beta+2\alpha+2\gamma & 0 & 0 & 0\\
	0 & 0 & 0 & \mu_l & 0 & 0\\
	0 & 0 & 0 & 0 & \mu_l & 0\\
	0 & 0 & 0 & 0 & 0 & \mu_t
\end{pmatrix}
\begin{pmatrix}
	\varepsilon_{11} \\ \varepsilon_{22} \\ \varepsilon_{33} \\ 2\varepsilon_{23} \\ 2\varepsilon_{13} \\ 2\varepsilon_{12}
\end{pmatrix}.
\ee
%==========================================================================================
%
%==========================================================================================
\subsection{The compliance tensor}
It is also useful to have available the inverse of \eqref{cauchy_tensor}, and to find the compliance tensor $\mathbb{S}$ corresponding to $\mathbb{C}$; that is,
\be\label{strain-tensor}
\bm{\varepsilon} = \mathbb{S} \bm{\sigma}.
\ee
We derive here an explicit form for $\mathbb{S}$, using an approach that is more direct than that in \cite{Lubarda-Chen2008}. First, from \eqref{cauchy_tensor}, 
\begin{align}\label{1st_varepsilon_sigma}
	\begin{pmatrix}
		\tr\bm{\varepsilon} \\ \bm{M}:\bm\varepsilon
	\end{pmatrix}
	= \dfrac{1}{\mathcal{K}}
	\begin{pmatrix}
		2\mu_t + \beta + \alpha + 2\gamma & -\beta - 3\alpha - 2\gamma\\
		-\lambda - \alpha & 3\lambda + 2\mu_t + \alpha
	\end{pmatrix}
	\begin{pmatrix}
		\tr{\bm\sigma} \\ \bm{M}:\bm{\sigma}
	\end{pmatrix},
\end{align}
and
\begin{align}\label{varepsilonM}
	\bm{\varepsilon M} = \dfrac{1}{2\mu_t + \gamma}(\bm{\sigma M} - \bm{M\sigma}) + \bm{M\varepsilon},
\end{align}
where
\be\label{K}
\mathcal{K} = 2(\lambda + \alpha)(\mu_t -\alpha) + 2(\lambda + \mu_t)(2\mu_t+\beta+\alpha+2\gamma).
\ee

Substituting back into \eqref{cauchy_tensor} and using $(\mu_t\bm{I} + \gamma \bm{M})^{-1} = \dfrac{1}{\mu_t(\mu_t +\gamma)} \big[(\mu_t+\gamma)\bm{I} - \gamma \bm{M}\big]$, equation \eqref{strain-tensor} becomes
\begin{align}\label{strain_stress}
\bm\varepsilon = \dfrac{1}{2\mu_t(\mu_t +\gamma)} &
	\Bigg[ (\mu_t+\gamma)\bm\sigma - \gamma \bm{M\sigma} - (\mu_t+\gamma)(A\lambda+C\alpha) (\tr\bm\sigma)\,\bm{I} \nonumber \\
& 	 - \left[\mu_t(A\alpha+C\beta) - \gamma (A\lambda+C\alpha)\right](\tr\bm\sigma)\,\bm{M}\nonumber\\
	& - (\mu_t+\gamma)(B\lambda+D\alpha) \bm{I}\bm{(M:\sigma)}
	 - \left[\mu_t(B\alpha+D\beta) - \gamma (B\lambda+D\alpha)\right]\bm{M}\bm{(M:\sigma)} \nonumber \\
	& -\dfrac{\gamma(\mu_t+\gamma)}{2\mu_t+\gamma} (\bm{\sigma M} - \bm{M\sigma})
	+ \dfrac{\gamma^2}{2\mu_t+\gamma} (\bm{M\sigma M} - \bm{M\sigma})\Bigg]\nonumber,
\end{align}
where
\be\label{ABCD}
	A = \dfrac{2\mu_t+\beta+\alpha+2\gamma}{\mathcal{K}},
	\hspace{0.5cm} B = -\dfrac{\beta + 3\alpha + 2\gamma}{\mathcal{K}},
	\hspace{0.5cm} C = -\dfrac{\lambda + \alpha}{\mathcal{K}}
	\hspace{0.5cm} \text{and} \hspace{0.5cm}
	D = \dfrac{3\lambda + 2\mu_t + \alpha}{\mathcal{K}}.
\ee
Next, we obtain expressions for the five material parameters in \eqref{cauchy_tensor} in terms of five physically meaningful constants, viz.
$E_t$: Young's modulus in the transverse direction; $E_l$: Young's modulus in the fibre direction; and
$\nu_t$ and $\nu_l$: respectively Poisson's ratios for the transverse strain with respect to the fibre direction and the plane normal to it.
The remaining constants are the two shear moduli $\mu_t$ and $\mu_l$, and one may further define  $\mu_t$ by 
\be
\mu_t = \dfrac{E_t}{2(1+\nu_t)}.
\label{GT}
\ee

Then, choosing the fibre direction to coincide with the basis vector $\bm{e}_3$, the compliance relation has the alternate form \cite{Exadaktylos2001}
\be\label{strain_stress_engineer_csts}
\begin{pmatrix}
	\varepsilon_{11} \\ \varepsilon_{22} \\ \varepsilon_{33} \\ 2\varepsilon_{23} \\ 2\varepsilon_{13} \\ 2\varepsilon_{12}
\end{pmatrix}
= 
\begin{pmatrix}
	\dfrac{1}{E_t} & -\dfrac{\nu_t}{E_t} & -\dfrac{\nu_l}{E_l} & 0 & 0 & 0\\
	-\dfrac{\nu_t}{E_t} & \dfrac{1}{E_t} & -\dfrac{\nu_l}{E_l} & 0 & 0 & 0\\
	-\dfrac{\nu_l}{E_l} & -\dfrac{\nu_l}{E_l} & \dfrac{1}{E_l} & 0 & 0 & 0\\
	0 & 0 & 0 & \dfrac{1}{\mu_l} & 0 & 0\\
	0 & 0 & 0 & 0 & \dfrac{1}{\mu_l} & 0\\
	0 & 0 & 0 & 0 & 0 & \dfrac{1}{\mu_t}
\end{pmatrix}
\begin{pmatrix}
	\sigma_{11} \\ \sigma_{22} \\ \sigma_{33} \\ \sigma_{23} \\ \sigma_{13} \\ \sigma_{12}
\end{pmatrix}.
\ee
The nature of the Poisson's ratios may be determined by considering some simple loading cases. For the case of uniaxial stress in the $x_3$-direction, with $\bm{a} = \bm{e}_3$, the strains are given by 
\be
\varepsilon_{11} =-\dfrac{\nu_l}{E_l}\sigma_{33}, \quad  \varepsilon_{22} =-\dfrac{\nu_l}{E_l}\sigma_{33} \quad \text{and} \quad \varepsilon_{33} = \dfrac{\sigma_{33}}{E_l}\,.
\label{nu_t}
\ee
Thus 
\be
\varepsilon_{11} = \varepsilon_{22} = -\nu_l\varepsilon_{33}\,,
\label{varepsilon1133}
\ee
so that $\nu_l$ determines the lateral contraction in the plane of isotropy, as a result of strain in the longitudinal or fibre direction.

Similarly, considering the case of uniaxial stress in the (transverse) $x_1$-direction, we have
\be
\varepsilon_{11} = \dfrac{\sigma_{11}}{E_t}, \quad \varepsilon_{22} = -\dfrac{\nu_t}{E_t}\sigma_{11}, \quad \text{and} \quad \varepsilon_{33} =-\dfrac{\nu_l}{E_l}\sigma_{11},
\label{nu_l}
\ee
so that 
\be
\varepsilon_{33} =-\nu_l\dfrac{E_t}{E_l}\varepsilon_{11}\,.
\label{varepsilon3311}
\ee
Thus, lateral contraction in the fibre direction depends directly on the ratio of Young's moduli in the fibre and transverse directions. In the inextensible limit, when $E_l/E_t \rightarrow \infty$, there is no lateral contraction in the fibre direction. 

By inverting \eqref{strain_stress_engineer_csts} and comparing with \eqref{stress_strain_e3}, the material parameters $\lambda,\,\alpha$ and $\beta$ can be written in terms of the engineering constants as
\begin{align}
	\lambda &= \dfrac{E_t(\nu_l^2 E_t + \nu_t E_l)}{(1+\nu_t)(E_l(1-\nu_t) - 2\nu_l^2 E_t)},\nonumber\\[10pt]
	\alpha  &= \dfrac{E_t\left[E_l \nu_l (1+\nu_t) - \nu_l^2 E_t - \nu_t E_l\right]}{(1+\nu_t)(E_l(1-\nu_t) - 2\nu_l^2 E_t)},\label{normalitoeng}\\[10pt]
	\beta   &= \dfrac{E_l^2(1 - \nu_t^2) - E_t^2\nu_l^2 + E_t E_l(1-2\nu_t\nu_l-2\nu_l)}{(1+\nu_t)(E_l(1-\nu_t) - 2\nu_l^2 E_t)} - 4\mu_l.\nonumber
\end{align}

Henceforth, we set
\be\label{cpq}
%\nu_l = c\nu_t , \hspace{0.3cm}
E_l = p E_t \hspace{0.3cm}\text{ and } \hspace{0.3cm} \mu_l = q \mu_t\,.
\ee
Thus $p$ measures the stiffness in the fibre direction relative to that in the plane of isotropy, and q is the ratio of the two shear moduli.
Using \eqref{cpq} the expressions \eqref{normalitoeng} become
\begin{subequations}\label{material_parameters}
\begin{align}
\mu_l &= \dfrac{qE_t}{2(1+\nu_t)},\\
\lambda &= \dfrac{(\nu_tp+\nu_l^2)}{(1+\nu_t)((1-\nu_t)p-2\nu_l^2)}E_t,\\[10pt]
\alpha &= \dfrac{(\nu_l-\nu_t+\nu_t\nu_l)p-\nu_l^2}{((1+\nu_t)((1-\nu_t)p-2\nu_l^2)}E_t,\\[10pt]
\beta &= \dfrac{(1-\nu_t^2)p^2 + (-2\nu_t\nu_l+2q\nu_t-2\nu_l+1-2q)p - (1-4q)\nu_l^2}{(1+\nu_t)((1-\nu_t)p-2\nu_l^2)}E_t.\label{beta}
\end{align}
\end{subequations}

\subsection{Pointwise stability}\label{Sec:pointwise_stability}
The condition of strong convexity, or pointwise stability, is equivalent to the positive definiteness of the strain energy \cite{Marsden-Hughes1994}: that is,
\be
\bm\varepsilon_{AB}\mathbb{C}_{ABCD}\bm\varepsilon_{CD}>0\quad\mbox{for any non-zero second-order tensor  $\bm\varepsilon$.}
\label{PS}
\ee

From \eqref{altC_ABCD} we have
\begin{align*}
\bm\varepsilon_{AB}\mathbb{C}_{ABCD}\bm\varepsilon_{CD}
&= \lambda (\tr \bm\varepsilon)^2 + 2\mu_t (\bm\varepsilon:\bm\varepsilon) + \beta (\bm\varepsilon:\bm{M})^2 + 2\alpha(\tr \bm\varepsilon)(\bm\varepsilon:\bm{M}) + 2\gamma(\bm\varepsilon:\bm{\varepsilon M})\\
&= \left(\lambda+\frac{2}{3}\mu_t\right) (\tr \bm\varepsilon)^2 + 2\mu_t |\mathrm{dev} \bm\varepsilon|^2 + \beta (\bm\varepsilon:\bm{M})^2\\
&\quad + 2\alpha(\tr \bm\varepsilon)(\bm\varepsilon:\bm{M}) + 2\gamma(\bm\varepsilon:\bm{\varepsilon M}).
\end{align*}
Here $\mbox{dev}\,\bm{\varepsilon} = \bm{\varepsilon} - (1/3)(\mbox{tr}\,\bm{\varepsilon})\bm{I}$ is the deviatoric part of the strain.
We assume that
\be
E_t,E_l>0, \text{ giving } p>0,
\ee
we further assume that
\be
\mu_l\geq\mu_t>0, \text{ implying } \gamma\geq 0
\ee
and from equation \eqref{GT}, $\nu_t>-1$.

We note that $\bm\varepsilon:\bm{M} = \bm\varepsilon\bm{a}\cdot\bm{a}$ and $\bm\varepsilon:\bm\varepsilon\bm{M} = |\bm\varepsilon\bm{a}|^2$, with
$|\bm\varepsilon\bm{a}|^2 \geq (\bm\varepsilon\bm{a}\cdot\bm{a})^2$,
then
\begin{align}\label{point_stab_eq}
\bm\varepsilon_{AB}\mathbb{C}_{ABCD}\bm\varepsilon_{CD} \geq&
	(\lambda+\mbox{$\frac{2}{3}$}\mu_t ) (\tr \bm\varepsilon)^2
	+ (\beta + 2\gamma) (\bm\varepsilon\bm{a}\cdot\bm{a})^2
	+ 2\alpha(\tr \bm\varepsilon)(\bm\varepsilon\bm{a}\cdot\bm{a})\nonumber\\
&	+ 2\mu_t|\mathrm{dev} \bm\varepsilon|^2\,.
\end{align}
We can obtain the pointwise stability condition by establishing conditions such that the right hand side of the inequality \eqref{point_stab_eq} is strictly positive; that is
\(
\mathcal{A} + 2\mu_t|\mathrm{dev} \bm\varepsilon|^2 > 0,
\)
with
\be
\mathcal{A} = \left(\lambda+\frac{2}{3}\mu_t\right) (\tr \bm\varepsilon)^2
	+ (\beta + 2\gamma) (\bm\varepsilon\bm{a}\cdot\bm{a})^2
	+ 2\alpha(\tr \bm\varepsilon)(\bm\varepsilon\bm{a}\cdot\bm{a}).
\label{calA}
\ee
The discriminant of $\mathcal{A}$, treated as a quadratic function of $\tr \bm\varepsilon$ and $\bm\varepsilon\bm{a}\cdot\bm{a}$, is
\be\label{delta}
\Delta = 4\left[\alpha^2 - \left(\lambda+\frac{2}{3}\mu_t\right)(\beta + 2\gamma)\right]
\ee
We assume that 
\be
\lambda+\frac{2}{3}\mu_t >0\,,
\label{lambdamupos}
\ee
so that if $\Delta <0$, then $\mathcal{A}>0$ and we have pointwise stability.

We have
\begin{equation}
\lambda + \dfrac{2}{3}\mu_t = \dfrac{(2\nu_t + 1)p + \nu_l^2}{3(1+\nu_t)[(1-\nu_t)p-2\nu_l^2]}E_t>0,
\end{equation}
which is true iff
\begin{equation}\label{equiv_lam_mu>0}
\begin{dcases*}
(2\nu_t + 1)p + \nu_l^2>0 \text{ and } (1-\nu_t)p-2\nu_l^2>0,\\
\text{OR }(2\nu_t + 1)p + \nu_l^2<0 \text{ and } (1-\nu_t)p-2\nu_l^2<0.
\end{dcases*}
\end{equation}
Now
\be\label{deno_cond}
(1-\nu_t)p-2\nu_l^2<0 \Longleftrightarrow \dfrac{(1-\nu_t)p}{2}<\nu_l^2.
\ee
If this holds and $(2\nu_t + 1)p + \nu_l^2<0$, it follows that
\be\label{num_cond}
\dfrac{(1-\nu_t)}{2} < -(2\nu_t+1) \Longleftrightarrow \nu_t+1<0,
\ee
which contradicts the assumption $\nu_t>-1$.
Therefore, \eqref{lambdamupos} is equivalent to the conditions
%\be
\begin{subequations}
\begin{align}
&(2\nu_t + 1)p + \nu_l^2>0,\label{lam_mu_cond_1}\\
\text{ and } &(1-\nu_t)p-2\nu_l^2>0.\label{lam_mu_cond_2}
\end{align}
\end{subequations}
%\ee
The second condition is the same as that obtained in \cite{Exadaktylos2001,Lai2009};
and it implies that $\nu_t<1$.

We return to \eqref{delta} for the discriminant, and by using the expressions of the material parameters in \eqref{material_parameters}, rewrite this as a function of $p$ as follows:
\begin{align*}
\Delta  &= \dfrac{4E_t^2}{3(1+\nu_t)^2\left((1 - \nu_t)p - 2\nu_l^2\right)^2} \Bigg[(\nu_t^2-1)\big( 2\nu_t + 1\big)p^3\\
		& +\Big\{\big[4\nu_t^2+6\nu_t+2\big]\nu_l^2 + \big[-2\nu_t^2+2\big]\nu_l -\nu_t^2 + 1\Big\}p^2
		+\Big\{\big[-4\nu_t - 4\big]\nu_l^3 + \big[- 2\nu_t -2\big]\nu_l^2\Big\}p \Bigg]\\
		&= -\dfrac{4E_t^2}{3(1+\nu_t)^2\left((1 - \nu_t)p - 2\nu_l^2\right)}
   			\Bigg( (\nu_t+1)\big(2\nu_t+1\big)p^2 - \Big[(2\nu_t+2)\nu_l + \nu_t + 1\Big]p  \Bigg)\\
		&= -\dfrac{4E_t^2p}{3(1+\nu_t)} \dfrac{(2\nu_t + 1)p - (2\nu_l + 1)}{(1 - \nu_t)p - 2\nu_l^2}.
\end{align*}

From this expression, with \eqref{lam_mu_cond_2}, we have $\Delta<0$  iff
\be\label{N}
(2\nu_t+1)p - (2\nu_l+1)>0.
\ee
Furthermore, \eqref{N} implies \eqref{lam_mu_cond_1} since
\[
0 \leq (\nu_l + 1)^2 = \nu_l^2+2\nu_l+1
\]
so that
\be\label{bound}
2\nu_l+1\geq -\nu_l^2.
\ee
Combining \eqref{bound} with \eqref{N}, we obtain \eqref{lam_mu_cond_1}.

We summarise as follows: sufficient conditions on the material constants for the elasticity tensor to be pointwise stable are
\begin{subequations}
\label{psconds}
\begin{align}
&\mu_l\geq\mu_t>0, \quad
p>0, \quad
\nu_t>-1,\\
& (2\nu_t+1)p - (2\nu_l+1)>0,\label{N_b}\\
\text{ and }& (1-\nu_t)p-2\nu_l^2>0.\label{D_b}
\end{align}
\end{subequations}

We illustrate these conditions with some examples.
\begin{enumerate}
\item[(a)] Assuming $\nu_l=\nu_t=\nu$, the zones of admissible values of $p$ and $\nu$ corresponding to the inequalities \eqref{N_b} and \eqref{D_b} are shown in the cross-hatched areas in Figure \ref{fig:admissible_values}.
\begin{figure}[h!]
\centering
\includegraphics[width=.5\columnwidth]{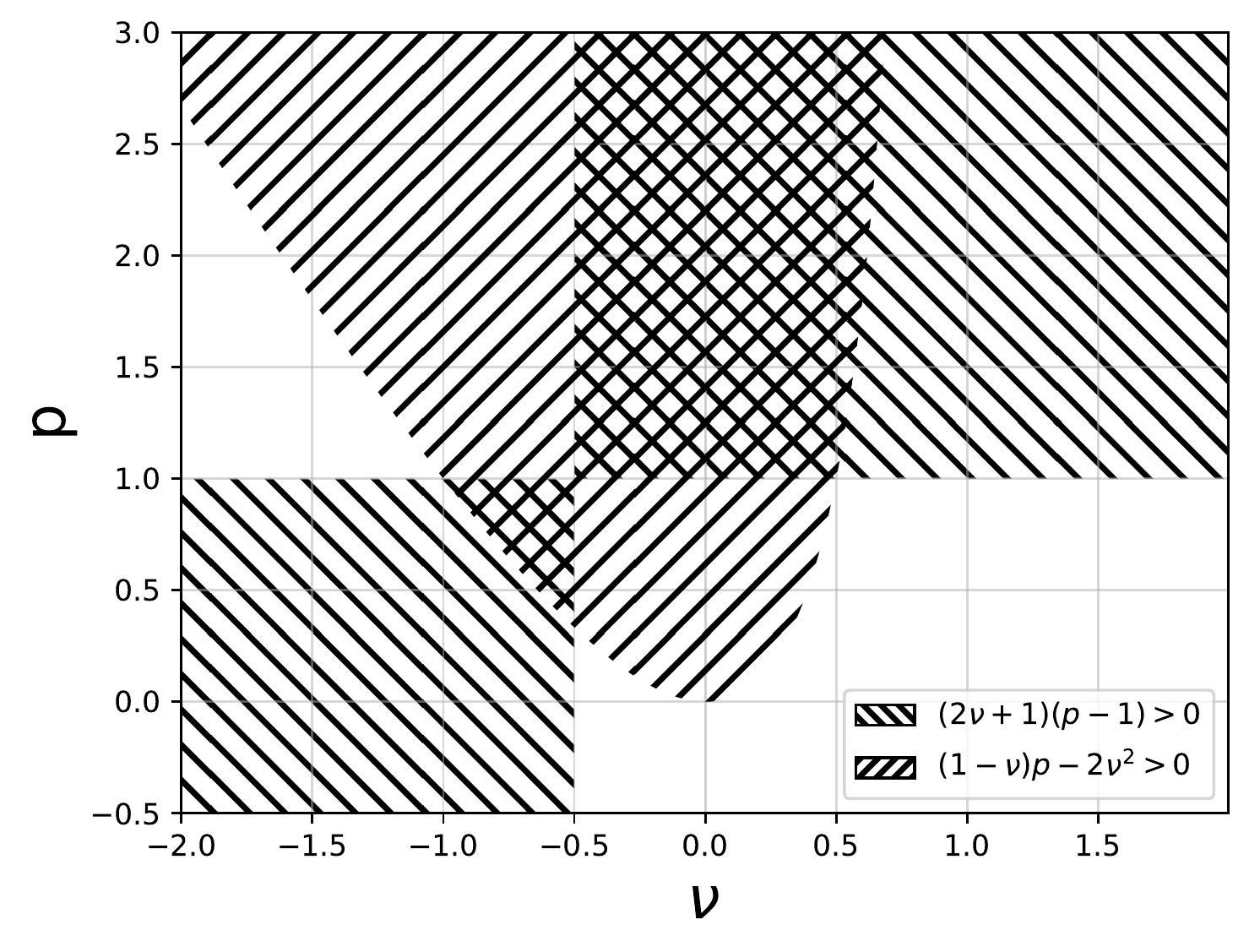}
\caption{Admissible values for $p$ and $\nu$, according to the inequalities \eqref{N_b} and \eqref{D_b}}
\label{fig:admissible_values}
\end{figure}

\item[(b)] If we choose $p=1$, then \eqref{N_b} is equivalent to
\[
2\nu_t+1> 2\nu_l+1 \Longleftrightarrow \nu_t > \nu_l
\]
Equation \eqref{D_b} is equivalent to
\[
1-\nu_t>2\nu_l^2 \Longleftrightarrow  |\nu_l| < \sqrt{\dfrac{1-\nu_t}{2}}
\]
Then, for example, for $\nu_t = 0.5$ we have $-0.5 < \nu_l < 0.5$.
\end{enumerate}

{\bf Remarks}
\begin{enumerate}
\item[1.]
The expression of $\Delta$ is independent of $q$; in fact,
\[
\beta+2\gamma = \dfrac{(1-\nu_t^2)p^2 + (-2\nu_t\nu_l+2\nu_t-2\nu_l-1)p + 3\nu_l^2}{(1+\nu_t)((1-\nu_t)p-2\nu_l^2)}E_t;
\]
thus conditions for pointwise stability hold for any value of $q\geq 1$.

\item[2.]
For the special case $\nu_t=\nu_l$ and $p=1$, we have $\beta+2\gamma = 0$ and $\alpha = 0$, so that
\[
\mathcal{A} = \left(\lambda+\frac{2}{3}\mu_t\right) (\tr \bm\varepsilon)^2,
\]
in which case $\lambda+\frac{2}{3}\mu_t>0$ is sufficient to ensure $\mathcal{A}>0$ provided only that $\tr\bm\varepsilon \neq 0$, and therefore $\lambda+\frac{2}{3}\mu_t>0$ and $\mu_l\geq\mu_t>0$ would ensure pointwise stability.
This includes the case of isotropy, for which $\mu_t = \mu_l$.
\end{enumerate}

%----------------------------------------------------------------------------------------
%	LINEAR_ELASTICITY
%----------------------------------------------------------------------------------------

\section{Governing equations and weak formulation}\label{Sec:Governing_eq}
Consider a transversely isotropic elastic body occupying a bounded domain $\Omega \subset \mathbb{R}^d,d=\{2,3\}$, with boundary $\Gamma = \Gamma_D \cup \Gamma_N$ having exterior unit normal $\bm{n}$. Here $\Gamma_D$ is the Dirichlet boundary, $\Gamma_N$ the Neumann boundary, and $\Gamma_D\cap \Gamma_N = \emptyset$.
The equilibrium equation is
\be\label{equilibrium}
	-\divg \bm\sigma(\bm{u}) = \bm{f}
\ee
and the boundary conditions are
\begin{subequations}\label{BC}
	\begin{align}
		\bm{u} = \bm{g} & \text{ on } \Gamma_D,\label{Dirichlet}\\
		\bm\sigma(\bm{u})\bm{n} = \bm{h} & \text{ on } \Gamma_N.\label{Neumann}
	\end{align}
\end{subequations}
Here $\bm\sigma$ is the Cauchy stress tensor defined by equation \eqref{cauchy_tensor}, $\bm{u}$ is the displacement vector, $\bm{f}$ is the body force, $\bm{g}$ a prescribed displacement, and $\bm{h}$ a prescribed surface traction.

We denote by $\mathcal{H}^1(\Omega)$ the Sobolev space of functions which, together with their generalized first derivatives, are square-integrable, and set
\[
\mathcal{V} = \{\bm{u} \in \big[\mathcal{H}^1(\Omega)\big]^d; \;\bm{u} = \bm{0} \; \text{on} \; \Gamma_D\},
\]
which is endowed with the norm
\[
\Vnorm{\cdot} = \|\cdot\|_{[\mathcal{H}^1(\Omega)]^d}.
\]
We will also use the norm
\[
\|\cdot\|_{[\mathcal{L}^2(\Omega)]^d} := \|\cdot\|_0\,.
%, \quad \|{\cdot}\|_{\mathcal{L}^2(\Omega_e)} = \|{\cdot}\|_{0,e}, \text{ and}\quad \quad \|{\cdot}\|_{\mathcal{L}^2(*)} = \|{\cdot}\|%_{0,*}.
\]
To take account of the non-homogeneous boundary condition \eqref{Dirichlet}, we define the function $\bm{u}_g \in [\mathcal{H}^1(\Omega)]^d$ such that $\bm{u}_g = \bm{g}$ on $\Gamma_D$, and the bilinear form $a(\cdot,\cdot)$ and linear functional $l(\cdot)$ by
\begin{subequations}
\begin{align}
a: [\mathcal{H}^1(\Omega)]^d \times [\mathcal{H}^1(\Omega)]^d \rightarrow \mathbb{R},& \hspace{1cm} a(\bm{u},\bm{v}) = \intO{\bm{\sigma}(\bm{u}):\bm{\varepsilon}(\bm{v})},\label{adef}\\
l: [\mathcal{H}^1(\Omega)]^d \rightarrow \mathbb{R},& \hspace{1cm} l(\bm{v}) = \intO{\bm{f}\cdot\bm{v}} + \intG{N}{\bm{h}\cdot\bm{v}} - a(\bm{u}_g,\bm{v}).\label{ldef}
\end{align}
\end{subequations}
The weak form of the problem is then as follows: given $\bm{f} \in [\mathcal{L}^2(\Omega)]^d$ and $\bm{h} \in [\mathcal{L}^2(\Gamma_N)]^d$,
find $\bm{U} \in [\mathcal{H}^1(\Omega)]^d$ such that $\bm{U}=\bm{u}+\bm{u}_g, \bm{u} \in \mathcal{V}$, and
\be\label{weak_form}
a(\bm{u},\bm{v}) = l(\bm{v}) \hspace{1cm} \forall \bm{v}\in \mathcal{V}.
\ee
We write the bilinear form as
\[
a(\bm{u},\bm{v}) = a^{iso}(\bm{u},\bm{v}) + a^{ti}(\bm{u},\bm{v}),
\]
where
\begin{subequations}\label{a_form}
\begin{align}
a^{iso}(\bm{u},\bm{v})
				&= \lambda \intO{(\nabla \cdot \bm{u})(\nabla \cdot \bm{v})}
					+ 2 \mu_t \intO{\bm{\varepsilon}(\bm{u}):\bm{\varepsilon}(\bm{v})},\label{a_iso}\\
a^{ti}(\bm{u},\bm{v})
					&=\alpha \intO{\left[(\bm{M}:\bm{\varepsilon}(\bm{u}))(\nabla \cdot \bm{v}) + (\nabla \cdot \bm{u})(\bm{M}:\bm{\varepsilon}(\bm{v}))\right]}
				+ \beta \intO{(\bm{M}:\bm{\varepsilon}(\bm{u}))(\bm{M}:\bm{\varepsilon}(\bm{v}))}\nonumber\\
					&+ \gamma \intO{\left[\bm{\varepsilon}(\bm{u})\bm{M}:\bm{\varepsilon}(\bm{v}) + \bm{M}\bm{\varepsilon}(\bm{u}):\bm{\varepsilon}(\bm{v})\right]}.\label{a_TI}
\end{align}
\end{subequations}
Note that $a^{iso}(\cdot,\cdot)$ and $a^{ti}(\cdot,\cdot)$ are symmetric.
The well-posedness of the weak problem requires the bilinear form to be continuous and coercive, and the linear functional continuous.

%=================================================================================
%
%=================================================================================
{\bf Assumption.} We assume the coefficients in the elasticity tensor $\mathbb{C}$ to satisfy the conditions for pointwise stability given by \eqref{psconds}.

\textbf{Continuity.}\hspace{1cm}
For all $\bm{u}, \bm{v} \in \mathcal{V}$, we have
\begin{align*}
|a^{iso}(\bm{u},\bm{v})| &= \left|\lambda \intO{(\nabla \cdot \bm{u})(\nabla \cdot \bm{v})} + 2 \mu_t \intO{\bm{\varepsilon}(\bm{u}):\bm{\varepsilon}(\bm{v})}\right|\\
						&\leq \max{(\lambda,2\mu_t)}\Vnorm{\bm{u}} \Vnorm{\bm{v}}.
\end{align*}
Next, we bound the first term on the right-hand side of \eqref{a_TI} as follows:
\begin{align*}
	 \alpha \left|\intO{(\bm{M}:\bm{\varepsilon}(\bm{u}))(\nabla \cdot \bm{v})}\right|
	&\leq \alpha \Znorm{\bm{M}:\bm{\varepsilon}(\bm{u})}  \Znorm{\nabla \cdot \bm{v}}\\
	&\leq \alpha C_\alpha \Znorm{\bm{\varepsilon}(\bm{u})}  \Znorm{\nabla \cdot \bm{v}}\\
	&\leq \alpha C_\alpha \Vnorm{\bm{u}} \Vnorm{\bm{v}}.
\end{align*}
The other terms are bounded similarly, and we find that $a$ is continuous; that is,
\be\label{continuity}
|a(\bm{u},\bm{v})| \leq M\Vnorm{\bm{u}} \Vnorm{\bm{v}},
\ee
where
\be\label{M}
M = C(\max{(\lambda,2\mu_t)+\alpha+\beta+\gamma}).
\ee
%===========================================================
%
%===========================================================
\textbf{Coercivity.}\hspace{1cm}
For all $\bm{v} \in \mathcal{V}$ and $0\leq k < 1$ we have
\begin{align}
& a(\bm{v},\bm{v}) = \intO{\bm\sigma(\bm{v}):\bm\varepsilon(\bm{v})}\nonumber\\
	&= \int_\Omega\Bigg[\left(\lambda+\frac{2k}{3}\mu_t\right) (\tr \bm\varepsilon)^2 + \mbox{$\frac{2}{3}$}(1-k)\mu_t (\tr \bm\varepsilon)^2 + 2\mu_t |\mathrm{dev} \bm\varepsilon|^2 + \beta (\bm\varepsilon:\bm{M})^2 + 2\alpha(\tr \bm\varepsilon)(\bm\varepsilon:\bm{M})\nonumber\\
	&\hspace{1cm} + 2\gamma(\bm\varepsilon:\bm{\varepsilon M})\Bigg] \;dx\nonumber\\
	&\geq \int_\Omega\Bigg[\left(\lambda+\frac{2}{3}k\mu_t\right) (\tr \bm\varepsilon)^2
	+ \mbox{$\frac{2}{3}$}(1-k)\mu_t (\tr \bm\varepsilon)^2
	+ (\beta+2\gamma) (\bm\varepsilon\bm{a}\cdot\bm{a})^2
	+ 2\alpha(\tr \bm\varepsilon)(\bm\varepsilon\bm{a}\cdot\bm{a}) \nonumber\\
	&\hspace{1cm} + 2\mu_t|\mathrm{dev} \bm\varepsilon|^2\Bigg]\ \;dx \nonumber\\
	&= \intO{\left[\mathcal{A}(k)+\dfrac{2}{3}(1-k)\mu_t (\tr \bm\varepsilon)^2
	+ 2\mu_t|\mathrm{dev} \bm\varepsilon|^2\right]}\label{av_v}
\end{align}	
where
\be
{\cal A}(k) = \left(\lambda+\frac{2}{3}k\mu_t\right) (\tr \bm\varepsilon)^2
	+ (\beta+2\gamma) (\bm\varepsilon\bm{a}\cdot\bm{a})^2
	+ 2\alpha(\tr \bm\varepsilon)(\bm\varepsilon\bm{a}\cdot\bm{a})\,.
\label{calAk}
\ee
We note that ${\cal A}(1) = {\cal A}$ as defined in \eqref{calA}, and thus ${\cal A}(1) > 0$ if the conditions \eqref{psconds} are satisfied. 
Since ${\cal A}(k)$ depends smoothly on $k$, it follows that ${\cal A}(1-\delta) > 0$ for sufficiently small $\delta > 0$. Choosing such a $\delta$, and $k = 1 - \delta$ in \eqref{av_v}, we obtain
\begin{align}
a(\bm{v},\bm{v})
	& \geq 2\mu_t \intO{\left(\delta |\mathrm{sph} \bm\varepsilon|^2 + |\mathrm{dev} \bm\varepsilon|^2\right)}\nonumber\\
	%& = \mathrm{min} \{{\blue 2}(1-k)\mu_t, 2\mu_t\} \left(\Lnorm{\mathrm{sph} \bm\varepsilon}^2 + \Lnorm{\mathrm{dev} \bm\varepsilon}^2\right)\\
	%& \geq \mathrm{min} \{{\blue 2}(1-k)\mu_t, 2\mu_t\} \;\Lnorm{\bm\varepsilon}^2 \hspace{2cm} \text{(Minkowski's inequality)}\\
	&\geq K \;\Vnorm{\bm{v}}^2, \label{coercivity}%\hspace{1cm} \text{(Korn's inequality)} 
\end{align}
in which $K = 2 \delta C \mu_t$, $C$ being a constant arising from Korn's inequality.
Hence $a$ is coercive.

Finally we have, using the trace theorem,
\begin{align*}
|l(\bm{v})| & \leq \Znorm{\bm{f}} \Znorm{\bm{v}} + \|\bm{h}\|_{0,\Gamma_N}\|\bm{v}\|_{0,\Gamma} + M \Vnorm{\bm{u}_g}\Vnorm{\bm{v}} \\
&\leq \Znorm{\bm{f}} \Vnorm{\bm{v}} + C_0\|\bm{h}\|_{0,\Gamma_N}\Vnorm{\bm{v}} + M \Vnorm{\bm{u}_g}\Vnorm{\bm{v}},\\
			&\leq C \Vnorm{\bm{v}}.
\end{align*}
Thus the problem \eqref{weak_form} has a unique solution.

%----------------------------------------------------------------------------------------
%----------------------------------------------------------------------------------------
%	CONFORMING FE APPROXIMATIONS
%----------------------------------------------------------------------------------------
\section{Conforming finite element approximations}
Suppose that $\Omega$ is polygonal (in $\mathbb{R}^2$) or polyhedral (in $\mathbb{R}^3$), and  partitioned into a shape-regular mesh comprising $n_e$ disjoint subdomains $\Omega_e$ with boundary $\partial\Omega_e$ and outward unit normal $\bm{n}_e$.
Denote by $\mathcal{T}_h := \{\Omega_e\}_e$ the set of all elements.
We define the discrete space $\mathcal{V}^h \subset \mathcal{V}$ by
\be\label{Vh_space}
\mathcal{V}^h = \{\bm{v}_h \in \mathcal{V} \cap \mathcal{C}(\bar{\Omega}) \;|\; \bm{v}_h|_{\Omega_e} \in [\mathcal{Q}_1(\Omega_e)]^d, \bm{v}_h =0 \text{ on } \Gamma_D\},
\ee
where ${\cal Q}_1(\Omega_e)$ is the space of polynomials on $\Omega_e$ of degree at most $1$ in each component.

The discrete problem corresponding to conforming approximations is as follows:
find $\bm{u}_h \in \mathcal{V}^h$ that satisfies
\be\label{CG_problem}
a(\bm{u}_h,\bm{v}_h) = l(\bm{v}_h) \hspace{1cm} \forall \bm{v}_h\in \mathcal{V}^h.
\ee
Since the bilinear form $a(\cdot,\cdot)$ is continuous and coercive and the linear functional $l(\cdot)$ is continuous, from standard finite element convergence theory \cite{Ciarlet1978} we have
\be\label{CG_error}
||\bm{u}-\bm{u}_h||_\mathcal{V} \leq C_1h,
\ee
in which, using \eqref{continuity} and \eqref{coercivity}, the constant $C_1$ is given by
\be\label{CG_constant}
C_1 = \dfrac{C(\max{(\lambda,2\mu_t)}+\alpha+\beta+\gamma)}{\mu_t}.
\ee

\begin{figure}[h!]
\centering
\includegraphics[width=.5\columnwidth]{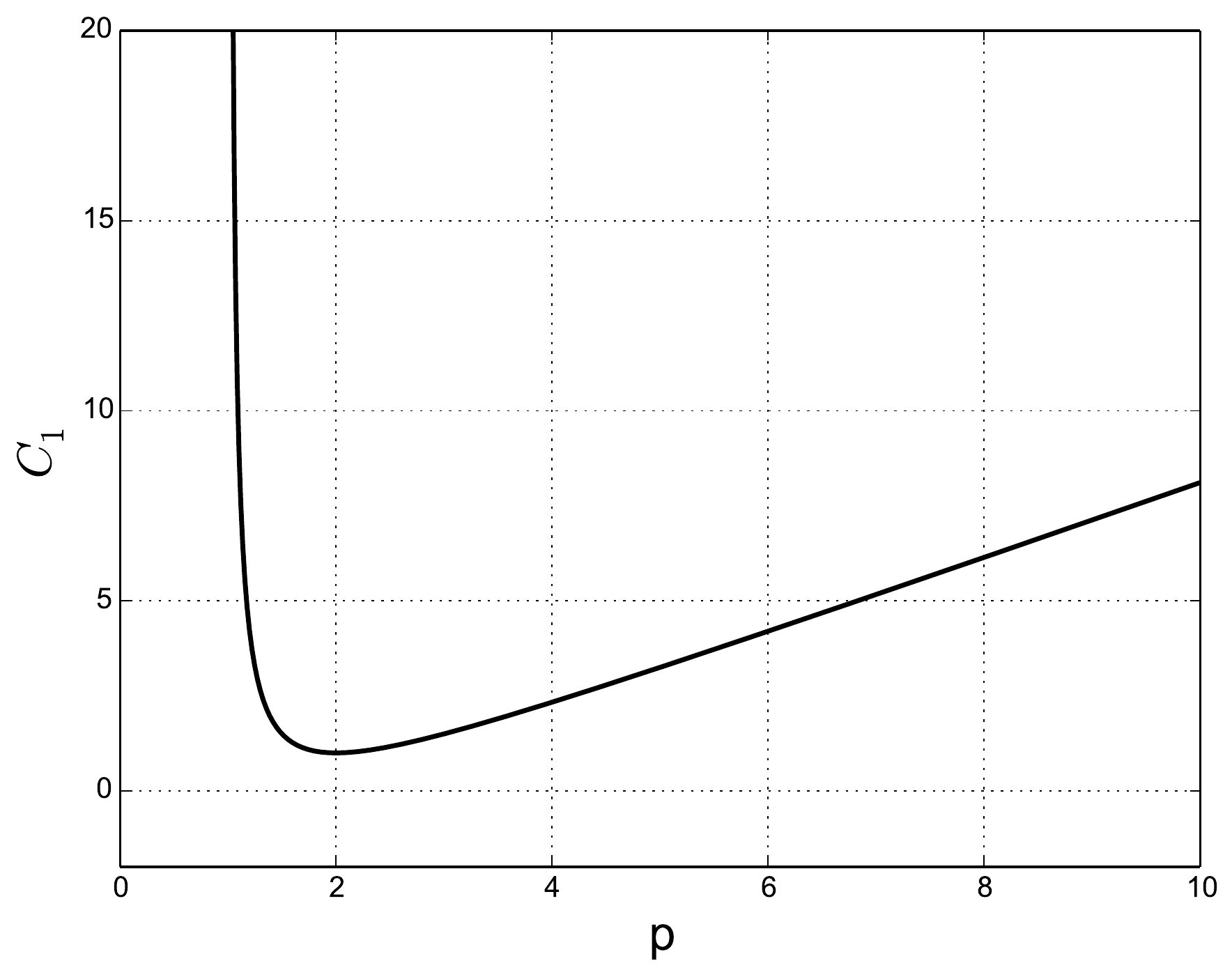}
\caption{The error bound constant $C_1$ in \eqref{CG_constant} against $p$, with $\nu_l=\nu_t=0.49995$ and $q =1$}
\label{fig:error_bound}
\end{figure}

From \eqref{material_parameters}, all three parameters $\lambda, \alpha,$ and $\beta$ have the same denominator $d := (1+\nu_t)((1-\nu_t)p-2\nu_l^2)$.
Fixing $\nu_l=\nu_t$ and $q =1$, for example, we have
\be\label{d}
d(p,\nu_t) = (1+\nu_t)((1-\nu_t)p-2\nu_t^2),
\ee
which tends to zero as
\[
p \rightarrow \frac{2\nu_t^2}{1 - \nu_t}.
\]
The limit $p=1$ corresponds to the case of isotropy, with the well-known limiting value $\nu_t = 1/2$, for which $\lambda$ becomes unbounded, while $\alpha = \beta = 0$.
For anisotropic materials, though, for which $p > 1$, it is seen from \eqref{CG_constant} and \eqref{d} that the constant $C_1$ in the error bound \eqref{CG_error} is bounded, so that in principle one has uniform convergence.
This is illustrated in Figure \ref{fig:error_bound}, which shows the behaviour of $C_1$ as a function of $p$, for near-incompressibility.
This issue will be explored later numerically, in Section \ref{sec:Numerical}.

\subsection{Under-integration}\label{Sec:UI}
The conforming approximation is known to exhibit volumetric locking while using $\mathcal{Q}_1$ elements (see for example \cite{Hughes1987}).
To overcome locking, we may under-integrate (i.e. use one-point integration) the terms involving volumetric and extensional deformation.
Let $\bar{\bm{x}}$ be the integration point and $\bar{\omega}$ the corresponding weight on element $\Omega_e$.

The volumetric term
\[
\lambda \int_{\Omega_e} (\nabla\cdot\bm{u}_h)(\nabla\cdot\bm{v}_h)\;dx
\]
is then replaced with
\[
\lambda (\nabla\cdot\bm{u}_h(\bar{\bm{x}}))(\nabla\cdot\bm{v}_h(\bar{\bm{x}})) \bar{\omega}.
\]
From \eqref{material_parameters}, we see that $\alpha$ is bounded as $p\rightarrow\infty$ (the inextensional limit), while $\beta \rightarrow \infty$ as $p\rightarrow\infty$.
Thus, in the event of extensional locking with $\mathcal{Q}_1$ elements, we will replace the extensional term
\[
\beta \int_{\Omega_e} (\bm{M}:\bm\varepsilon(\bm{u}_h))(\bm{M}:\bm\varepsilon(\bm{v}_h))\;dx
\]
with
\[
\beta \Big(\bm{M}:\bm\varepsilon(\bm{u}_h(\bar{\bm{x}}))\Big)\Big(\bm{M}:\bm\varepsilon(\bm{v}_h(\bar{\bm{x}}))\Big) \bar{\omega}.
\]

%\subsubsection*{Equivalence with mixed formulation}
One can easily show that one-point integration is the same as projection of the integrands onto the space of constants.
If we define by $\Pi_0$ the $\mathcal{L}^2$-orthogonal projection onto constants, then under-integrating the volumetric term is the same as replacing it with
\be\label{proj_vol}
\lambda \int_{\Omega_e} \Pi_0(\nabla\cdot\bm{u}_h)\Pi_0(\nabla\cdot\bm{v}_h)\;dx.
\ee
Similarly, the extensional term is replaced with
\be\label{proj_ext}
\beta \int_{\Omega_e} \Pi_0(\bm{M}:\bm\varepsilon(\bm{u}_h))\Pi_0(\bm{M}:\bm\varepsilon(\bm{v}_h))\;dx.
\ee
%Then we obtain the mixed displacement-pressure or $Q_1-P_0$ formulation.

\subsubsection*{Equivalence with perturbed Lagrangian formulation}
The perturbed Lagrangian approach takes as a starting point a strain energy of the form
  \be
  W (\bm{\varepsilon},T,\beta) = W^{\rm iso} + T(\bm{\varepsilon:M}) - \frac{1}{2\beta}T^2
  \label{PL}
  \ee
in which $W^{\rm iso}$ denotes that part of the strain energy that leads to \eqref{a_iso}, $\beta$ now plays the role of a penalty parameter, and $T$ is a Lagrange multiplier. For convenience we assume that $\alpha = 0$ in this section. The stress is obtained from
  \begin{align}
  \bm{\sigma} & = \frac{\partial W}{\partial \bm{\varepsilon}} \nonumber \\
& = \bm{\sigma}^{\rm iso} + T\bm{M} \label{stressPL}
 \end{align}
and in addition we have the condition
 \begin{align}
0 &= \frac{\partial W}{\partial T} \nonumber\\
& = \bm{\varepsilon:M} - \frac{T}{\beta}.\label{pressure}
 \end{align}

We show the equivalence of the perturbed Lagrangian formulation, used for example in \cite{Auricchio-Scalet-Wriggers2017}, to the under-integration formulation \eqref{proj_ext}.
In \eqref{cauchy_tensor} then for $T \in \mathcal{L}^2(\Omega)$ the extensional term in the weak formulation (see \eqref{a_TI}) is
\be\label{ext_stress}
\int_\Omega T (\bm{M}:\bm\varepsilon(\bm{v}))\,dx.
\ee
With a test function $\vartheta \in \mathcal{L}^2(\Omega)$, we can write the weak form of \eqref{pressure} as
\[
\int_{\Omega} \vartheta \Big(T - \beta (\bm{M}:\bm\varepsilon(\bm{u}))\Big) \,dx = 0.
\]
The discrete form with $\bm{u}_h \in \mathcal{Q}_1(\Omega_e)$ and $T_h,\,\vartheta_h \in \mathcal{P}_0(\Omega_e)$ at element level is
\[
\int_{\Omega_e} \vartheta_h \left( T_h - \beta (\bm{M}:\bm\varepsilon(\bm{u}_h))\right) \,dx = 0,
\]
and since $T_h$ is constant, we have
\[
T_h = \dfrac{\beta}{|\Omega_e|} \int_{\Omega_e} (\bm{M}:\bm\varepsilon(\bm{u}_h)) \,dx.
\]
We substitute in the discrete form of \eqref{ext_stress} to obtain
\begin{align*}\label{vol_PL}
\int_\Omega T_h (\bm{M}:\bm\varepsilon(\bm{v}_h))\,dx
&= \sum_{\Omega\in\mathcal{T}_h}\dfrac{\beta}{|\Omega_e|} \int_{\Omega_e} (\bm{M}:\bm\varepsilon(\bm{u}_h))\int_{\Omega_e} (\bm{M}:\bm\varepsilon(\bm{v}_h)) \;dx\\
&= \dfrac{\beta}{|\Omega_e|} \int_{\Omega_e} \Pi_0(\bm{M}:\bm\varepsilon(\bm{u}_h))\int_{\Omega_e} \Pi_0(\bm{M}:\bm\varepsilon(\bm{v}_h)) \;dx.
\end{align*}
Since $\bm{M}:\bm\varepsilon(\bm{u}_h)\in \mathcal{P}_1(\Omega_e)$, and is integrated exactly using one-point quadrature,
this expression is the same as \eqref{proj_ext}, so that the perturbed Lagrangian approximation is equivalent to a mixed $\mathcal{Q}_1-\mathcal{P}_0$ $(\bm{u}_h-T_h) $ formulation, and to under-integration of the extensional term.
%----------------------------------------------------------------------------------------
%	RESULTS
%----------------------------------------------------------------------------------------
\section{Numerical tests}\label{sec:Numerical}
In this section, we present the results of numerical simulations of two model problems to illustrate the formulations discussed in the preceding sections.
All examples are under conditions of plane strain and based on four- and nine-noded quadrilateral elements with standard bilinear and biquadratic interpolation of the displacement field.
%In particular, we are interested in demonstrating the locking behaviour that is consistent with the error bound.
Unless otherwise stated, all examples are presented for the case of near-incompressibility, i.e. the value of the transversal Poisson's ratio $\nu_t$ is close to $0.5$.
Precisely, we choose $\nu_t = 0.49995$.
We fix the value of the two Poisson's ratios to be equal, and also set $\mu_l = \mu_t$.
We consider values of $p>1$, so that the conditions \eqref{psconds} for pointwise stability are satisfied.

Define $\hat{a} := \widehat{(Ox,\bm{a})}$, the angle between the $x$-axis and the fibre direction $\bm{a}$.
For each problem, the following range of values for $\hat{a}$ will be considered:
\[
\hat{a} = \left\{0, \dfrac{\pi}{8}, \dfrac{\pi}{6}, \dfrac{\pi}{4}, \dfrac{\pi}{3}, \dfrac{3\pi}{8}, \dfrac{\pi}{2}, \dfrac{5\pi}{8}, \dfrac{3\pi}{4}, \dfrac{7\pi}{8}, \pi\right\}.
\]
Recall from Section \ref{Sec:UI} that under-integration is equivalent to a perturbed Lagrangian method.
Also, the standard formulation with large $\beta$ is equivalent to the penalty formulation as in \cite{Auricchio-Scalet-Wriggers2017}.
The results in the examples that follow are for the following element choices:
\begin{table}[h!]
\begin{tabular}{lp{12cm}}
$Q_1\_CG$ & The standard displacement formulation of order 1\\
$Q_2\_CG$ & The standard displacement formulation of order 2\\
$Q_1\_CG\_UI_\lambda$ & The standard displacement formulation with under-integration of the volumetric ($\lambda$-) term\\
$Q_1\_CG\_UI_\beta$ & The standard displacement formulation with under-integration of the extensional ($\beta$-) term\\
$Q_1\_CG\_UI_{\beta\lambda}$ & The standard displacement formulation with under-integration of the volumetric and the extensional terms
\end{tabular}
\end{table}

%====================================================================
%
%====================================================================
\subsection{Cook's membrane}
\begin{figure}[H]
\centering
\includegraphics[trim={0 5cm 0 4cm},clip,width=.35\columnwidth]{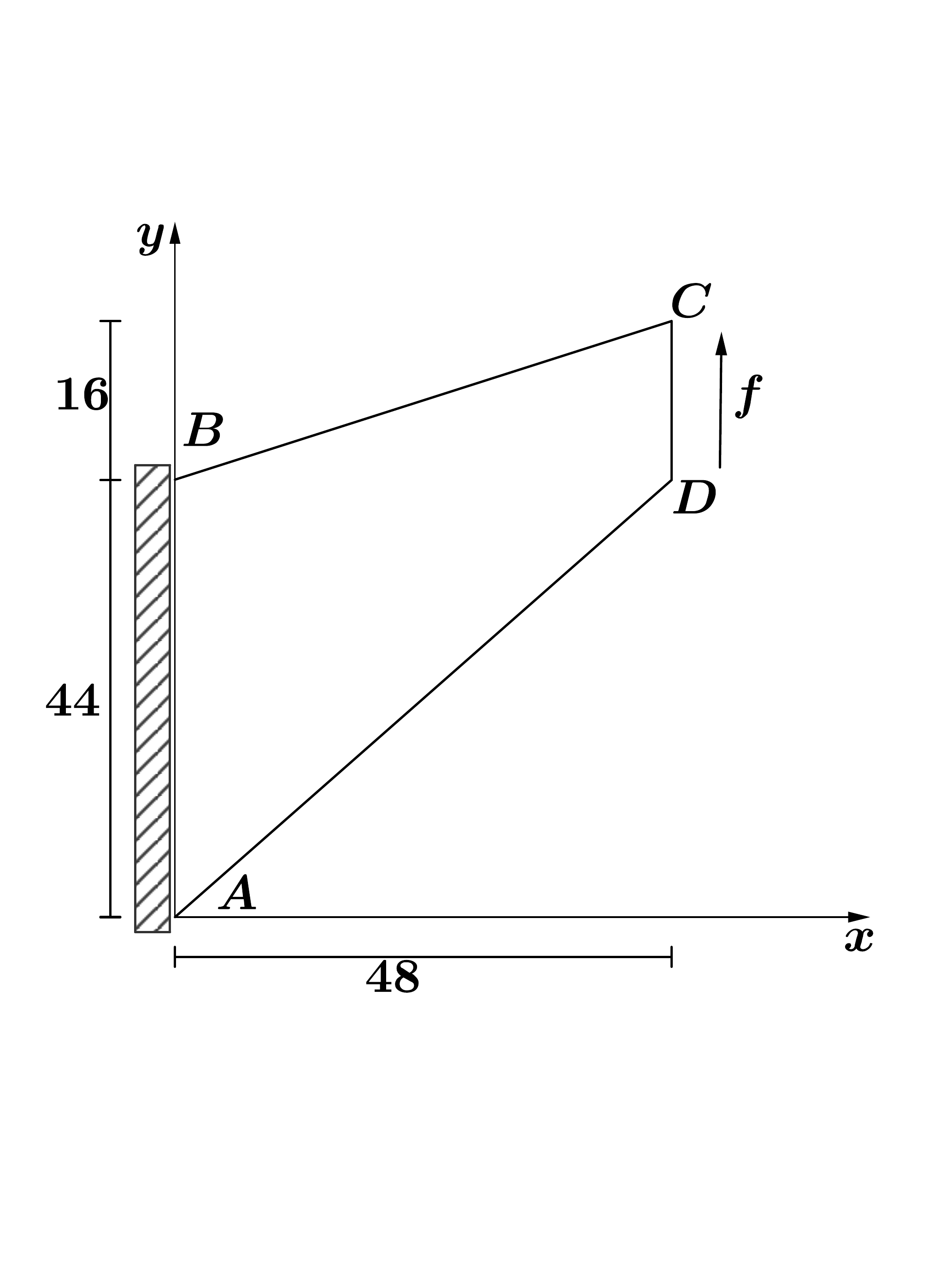}
\caption{Cook's membrane geometry and boundary conditions}
\label{fig:cooks}
\end{figure}
The Cook's membrane test consists of a tapered panel fixed on one edge and subject to a shearing load at the opposite edge as depicted in Figure \ref{fig:cooks}.
The applied load is $f=100$ and $E_t = 250$.
This test problem has no analytical solution.
The vertical tip displacement at corner $C$ is measured.
\begin{figure}[h!]
\centering
\subfloat[Moderate values of $p$]{\includegraphics[width=.7\columnwidth]{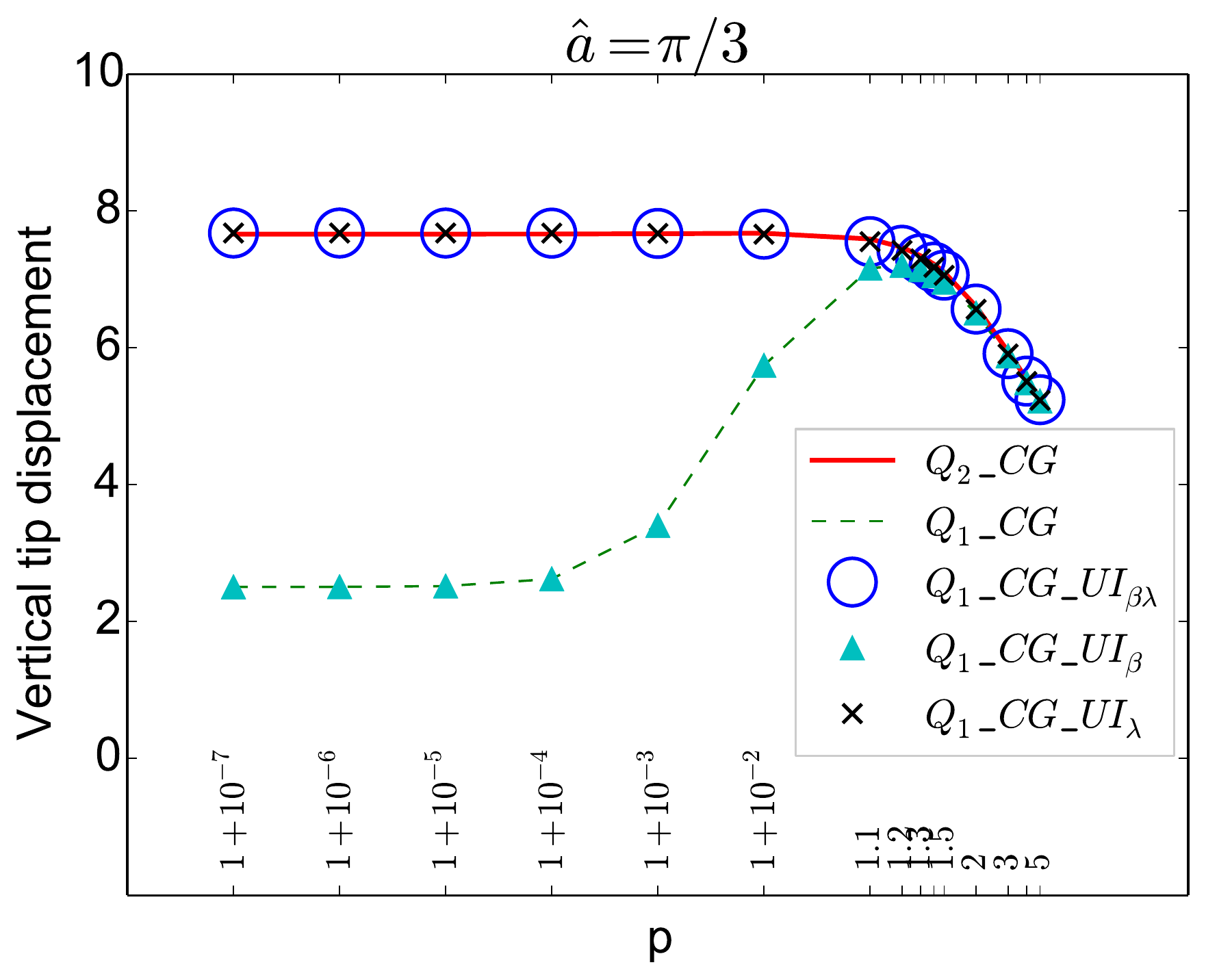}\label{fig:CM_pi3}}\\
\subfloat[High values of $p$]{\includegraphics[width=.7\columnwidth]{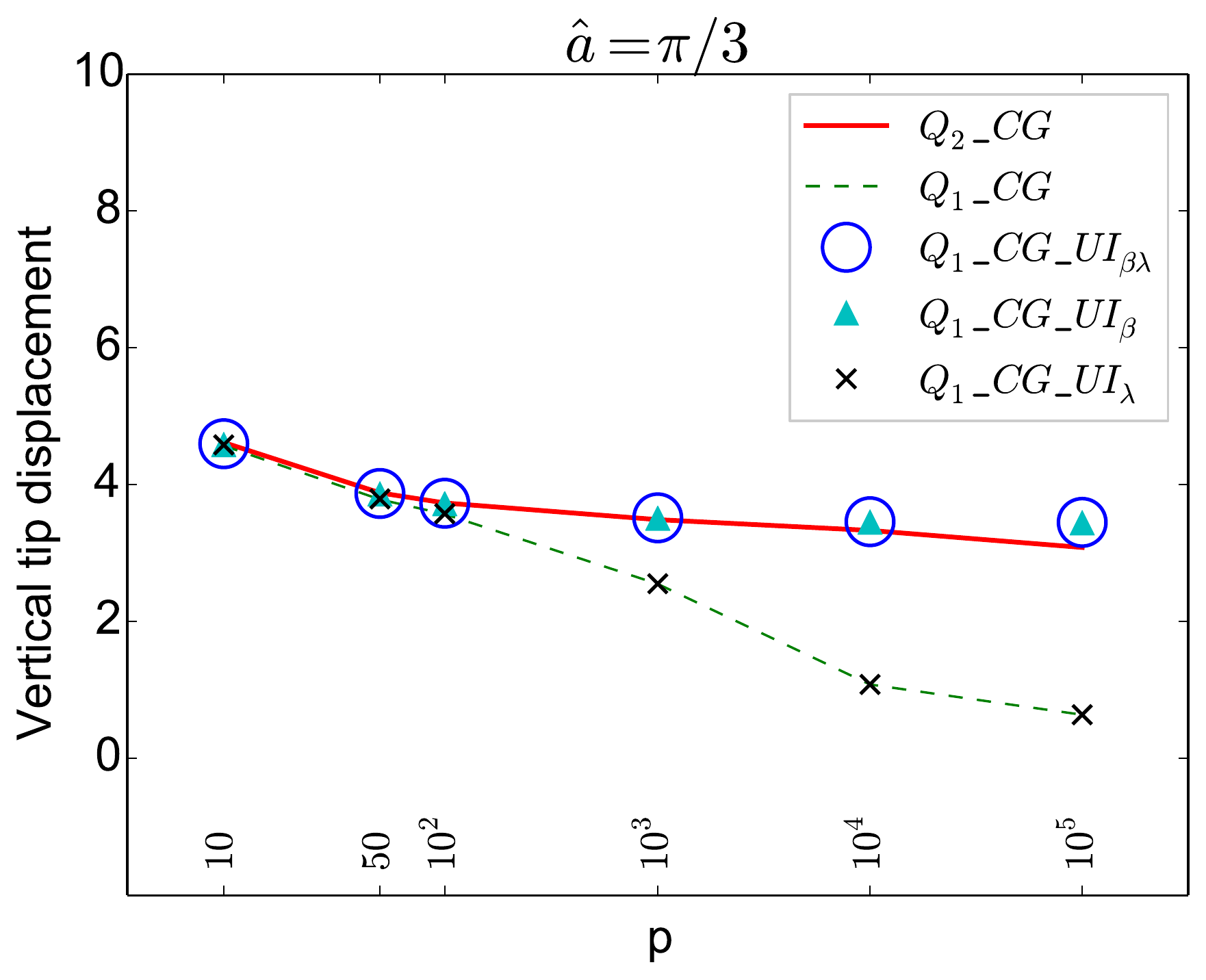}\label{fig:CM_big_pi3}}
\caption{Tip displacement vs $p$ for Cook's membrane problem.}
\label{fig:CM_tip_disp_log}
\end{figure}
Figure \ref{fig:CM_tip_disp_log} shows semilog plots of the tip displacement vs $p$ when the fibre direction is at the angle $\pi/3$ for the various element choices.
To investigate locking of the proposed formulation, we compare the results with the results obtained using the standard $Q_2$-element.

In Figure \ref{fig:CM_pi3}, for moderate values of $p\; (1\leq p \leq 5)$ the $Q_1\_CG$ formulation behaves well away from $p=1$, with evidence of locking behaviour as $p\rightarrow 1$.
Locking is avoided when the volumetric term is under-integrated $(Q_1\_CG\_UI_\lambda)$.
Notice that for $Q_1\_CG\_UI_\beta$ there is still locking, as the locking is purely volumetric.
In Figure \ref{fig:CM_big_pi3}, for higher values of $p\; (10\leq p \leq 10^5)$, the $Q_1\_CG$ method shows locking behaviour as $p$ get bigger, and convergent behaviour with $Q_1\_CG\_UI_\beta$.
Notice that under-integrating only the volumetric term $(Q_1\_CG\_UI_\lambda)$ has no effect since the locking is purely extensional.
$Q_1\_CG\_UI_{\beta\lambda}$ shows locking-free behaviour everywhere.
\begin{figure}
\centering
\includegraphics[width=.5\columnwidth]{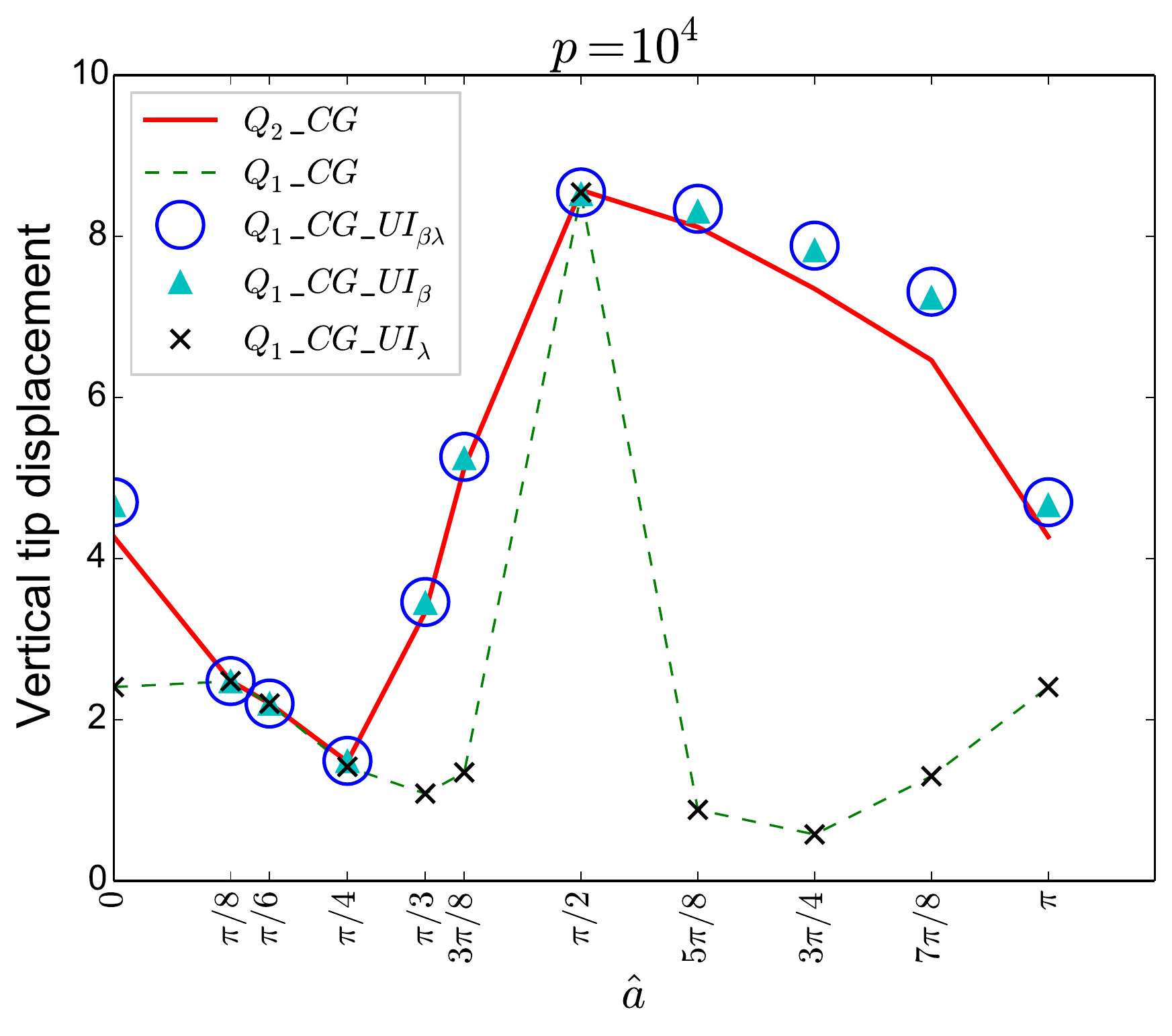}
\caption{Tip displacement measured at different fibre orientations, for $p=10^4$}
\label{fig:CM_orientation}
\end{figure}
Figure \ref{fig:CM_orientation} shows tip displacements for various fibre orientations, where the degree of anisotropy is fixed at $p=10^4$.
The behaviour observed in Figure \ref{fig:CM_tip_disp_log} is evident here over the range of fibre directions.
%====================================================================
%
%====================================================================
\subsection{Bending of a beam}
\begin{figure}[H]
\centering
\includegraphics[trim={0 10cm 0 10cm},clip,width=.5\columnwidth]{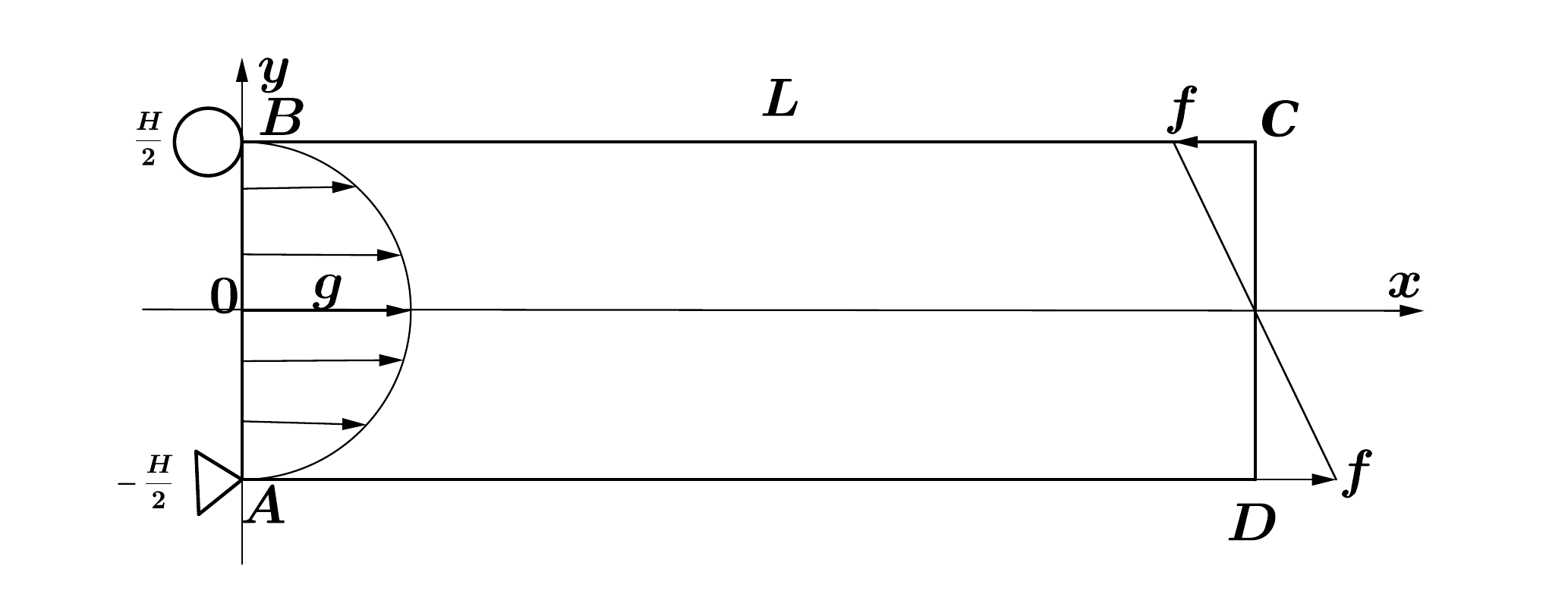}
\caption{Beam geometry and boundary conditions}
\label{fig:bending}
\end{figure}
We consider the beam shown in Figure \ref{fig:bending}, subject to a linearly varying load along the edge CD.
The horizontal displacement $\bm{u}$ is constrained at node B, while the node A is constrained in both directions.
The beam has length $L = 10$ and height $H = 2$ and the applied load has a maximum value $f = 3000$.
Here, $E_t = 1500$.
The boundary conditions are
\[
\begin{dcases}
u(0,y) = g(y),\\
v(0,-\frac{H}{2}) = 0,
\end{dcases}
\]
where
\[
g(y) = -\dfrac{f}{H}\mathbb{S}_{31} \left(y^2 - \dfrac{H^2}{4}\right).
\]
The compliance coefficients $\mathbb{S}_{ij}$ are given in the Appendix.
The analytical solution is
\[
\begin{dcases}
u(x,y) =& -\dfrac{2f}{H}\left( \mathbb{S}_{11}xy + \dfrac{1}{2}\mathbb{S}_{31} \left(y^2 - \dfrac{H^2}{4}\right) \right),\\
v(x,y) =& -\dfrac{f}{H}\left(\mathbb{S}_{21}\left(y^2 - \dfrac{H^2}{4}\right) - \mathbb{S}_{11}x^2 \right).
\end{dcases}
\]
The linearly varying load $\hat{f}$ with maximum $f$ is
\[
\hat{f}(y) = -\dfrac{2f}{H}y.
\]
\begin{figure}
\centering
\subfloat[Moderate values of $p$]{\includegraphics[width=.7\columnwidth]{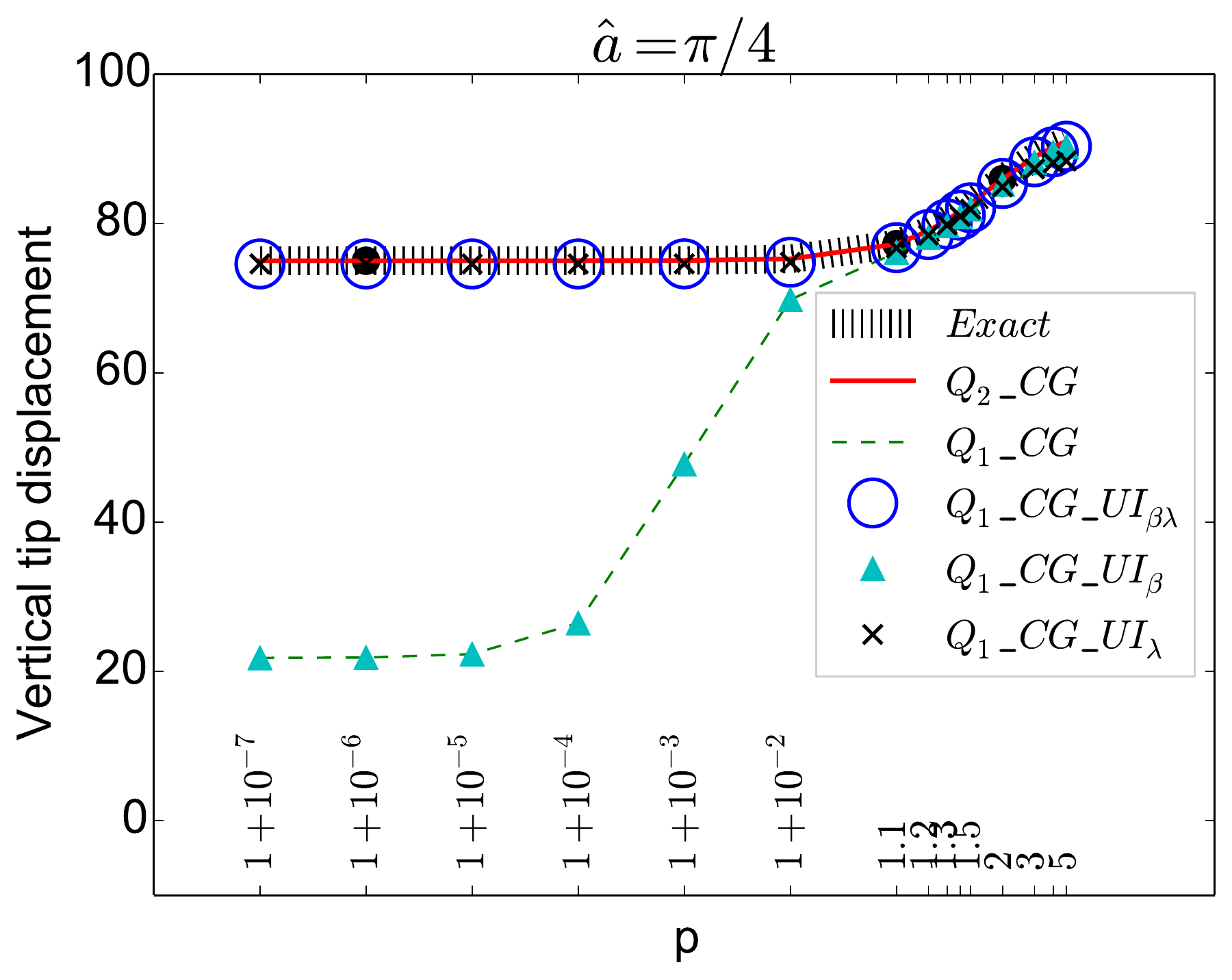}\label{fig:AB_pi4}}\\
\subfloat[High values of $p$]{\includegraphics[width=.7\columnwidth]{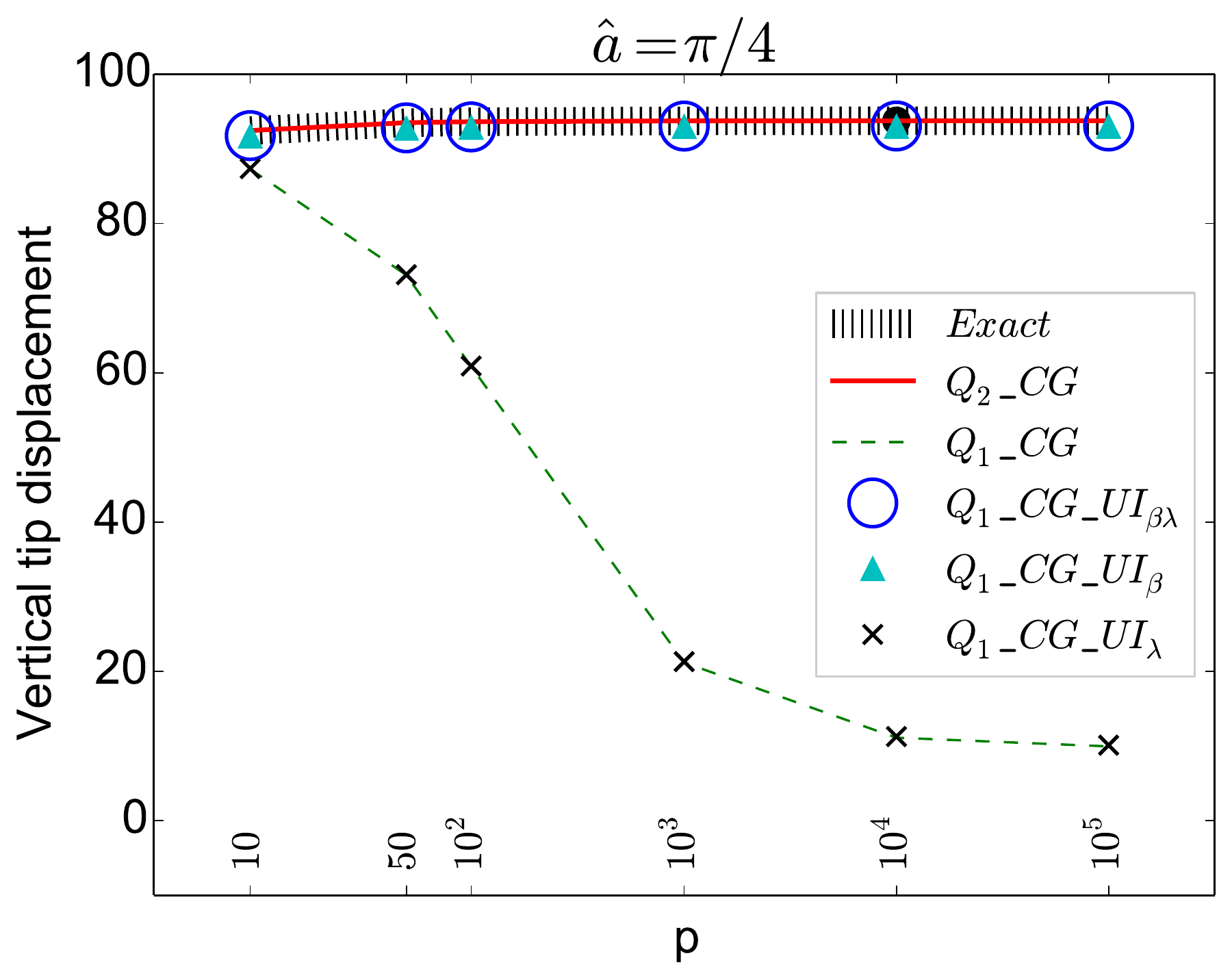}\label{fig:AB_big_pi4}}
\caption{Tip displacement vs $p$ for the beam problem} % The text in the square bracket is the caption for the list of figures while the text in the curly brackets is the figure caption
\label{fig:AB_tip_disp_log}
\end{figure}
In Figure \ref{fig:AB_tip_disp_log}, which shows semilog plots of tip displacement for different values of $p$, with the angle of the fibre direction $\pi/4$, locking behaviour is investigated by comparison with the analytical solution.
The same behaviour as appears for the Cook's example is seen, i.e. for moderate values of $p$ away from $p=1$ (approximately, $1.1 \leq p \leq 5$),
there is locking-free behaviour with $Q_1\_CG$, while locking occurs as $p$ approaches $1$. This is overcome by using under-integration $(Q_1\_CG\_UI_\lambda)$ (Figure \ref{fig:AB_pi4}).
For high values of $p$ there is purely extensional locking with $Q_1\_CG$, which is overcome by using $Q_1\_CG\_UI_\beta$ (Figure \ref{fig:AB_big_pi4}).
$Q_1\_CG\_UI_{\beta\lambda}$ shows locking-free behaviour for any value of $p$.
\begin{figure}
\centering
\includegraphics[width=.5\columnwidth]{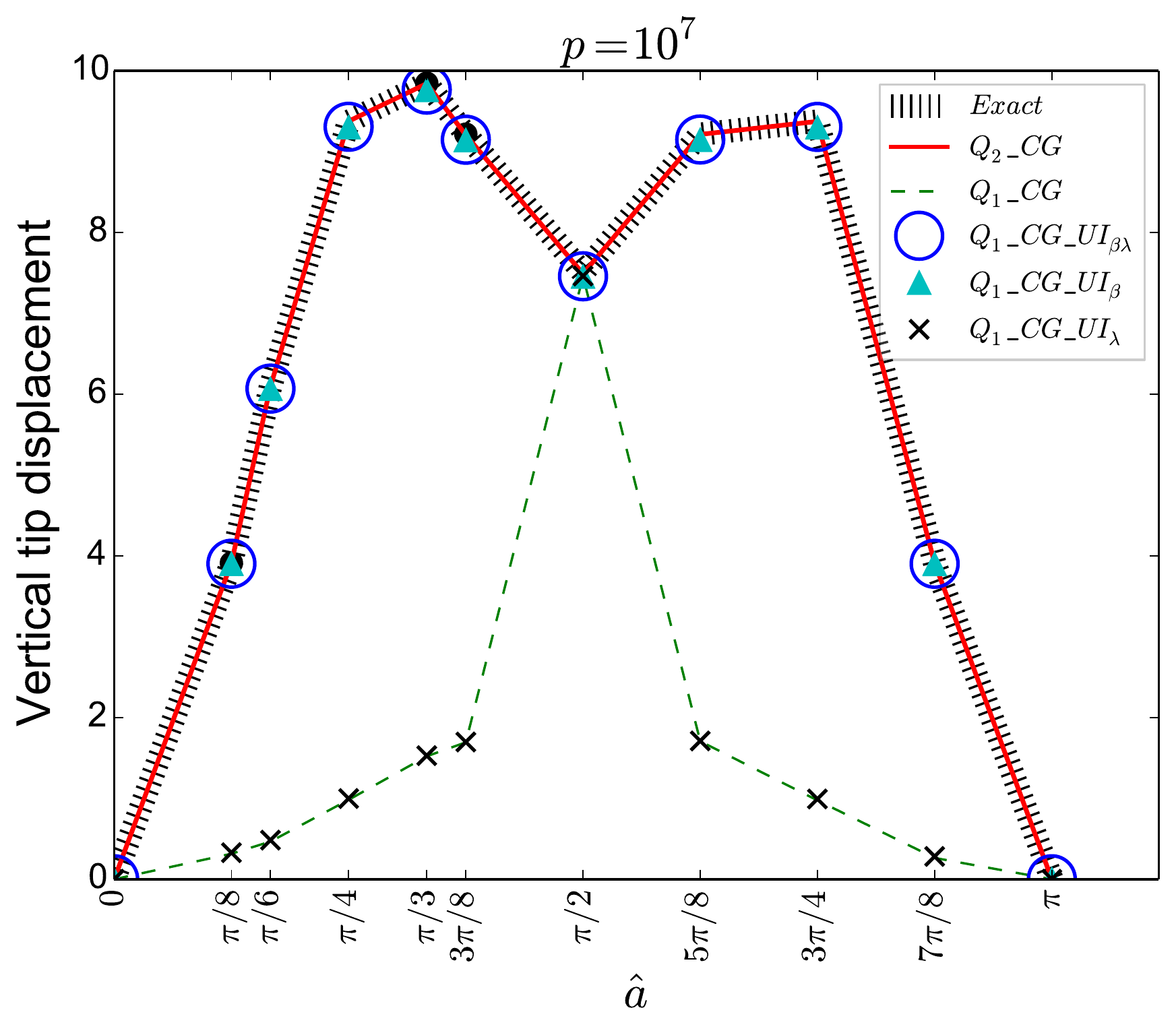}
\caption{Tip displacement measured at different fibre orientations}
\label{fig:AB_orientation}
\end{figure}
Figure \ref{fig:AB_orientation} shows tip displacements for various fibre orientations where the degree of anisotropy is fixed at $p=10^7$.
Extensional locking of $Q_1\_CG$ is observed except for the angles $0,\pi/2$.
For the angle $0\,(=\pi)$, the material is very stiff, relative to the type of loading.
%For the angle $\pi/2$ no locking would be expected as the extensional term of the formulation is
%\[
%\beta (\bm{M}:\bm\varepsilon) = \beta \varepsilon_{22} = \dfrac{\beta}{E_l}\big(\sigma_{22}-\nu_l\sigma_{11}\big),
%\]
%which is bounded when $p \rightarrow \infty$.

For the angle $\pi/2$ no locking is observed.
This can be accounted for by two factors: first, with this orientation the property of near-inextensibility in the vertical direction has a negligible effect on the bending-dominated deformation; and secondly, as previously discussed, the presence of anisotropy serves to circumvent volumetric locking.
\begin{figure}
\centering
\subfloat[$\hat{a} = \pi/4$]{\includegraphics[width=.5\columnwidth]{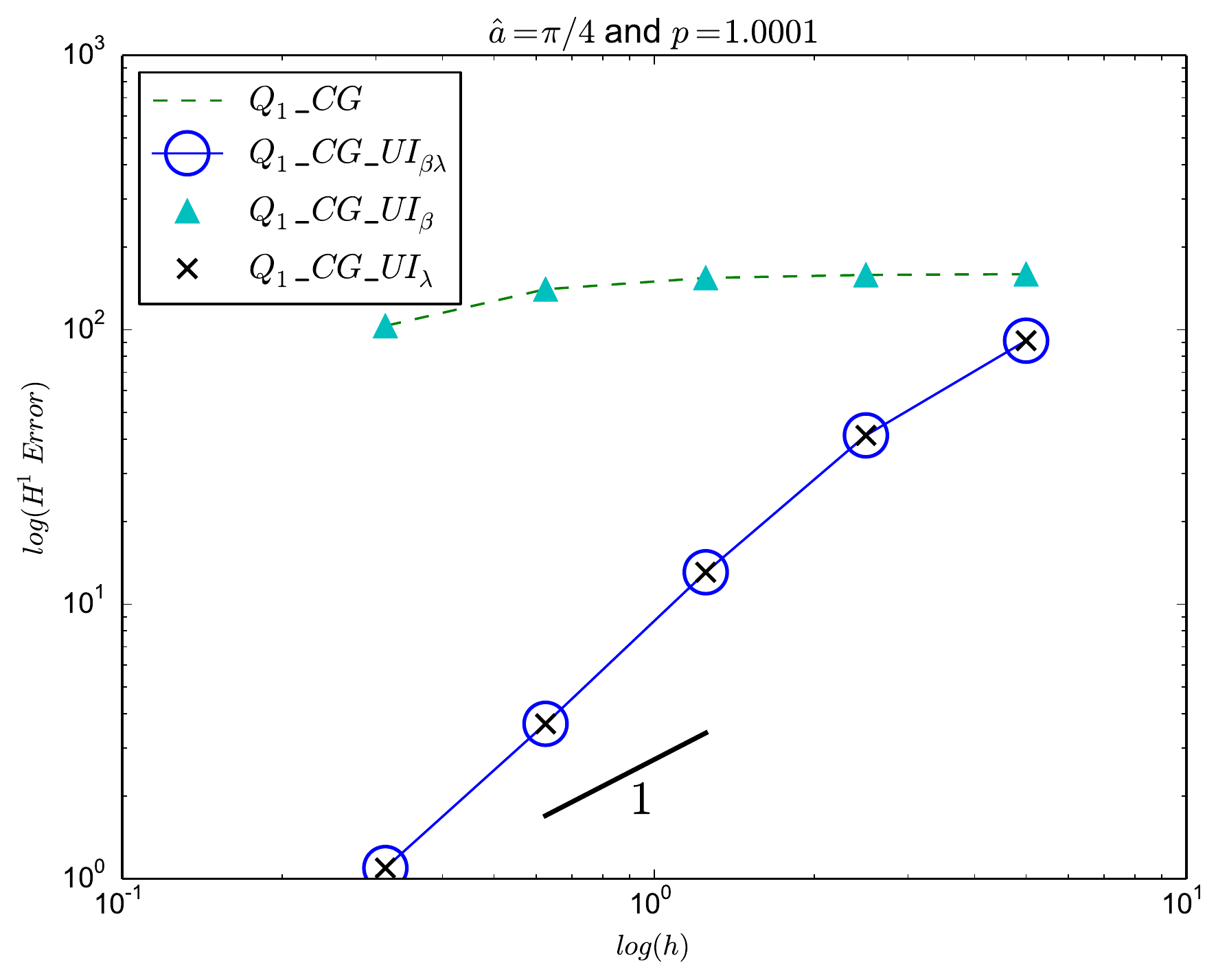}\label{fig:H1_p1_pi4}}
\subfloat[$\hat{a} = \pi/2$]{\includegraphics[width=.5\columnwidth]{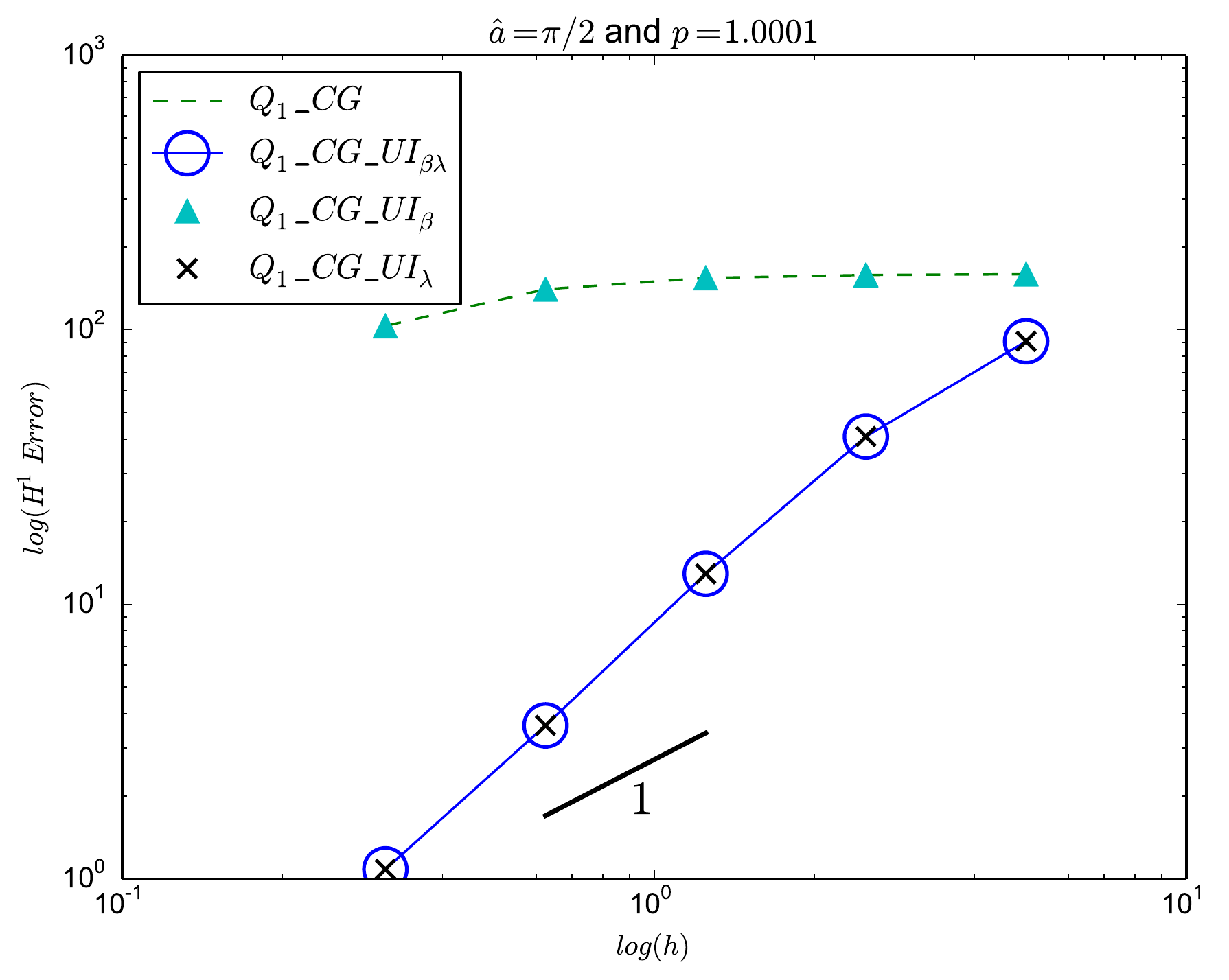}\label{fig:H1_p1_pi2}}\\
\subfloat[$\hat{a} = 3\pi/4$]{\includegraphics[width=.5\columnwidth]{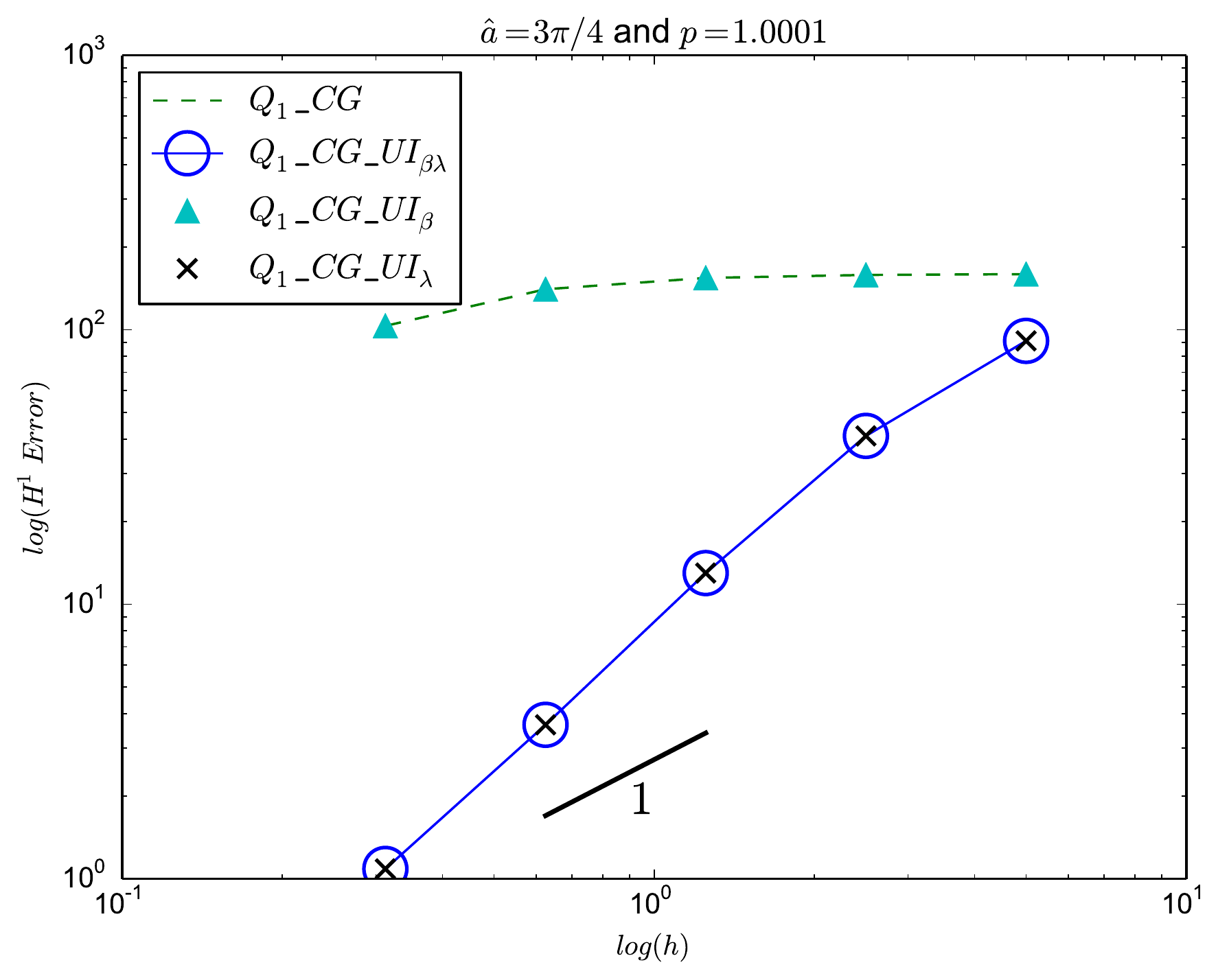}\label{fig:H1_p1_3pi4}}
\caption{Comparison of $\mathcal{H}^1$ errors for conforming and under-integrated elements on quadrilaterals, for different fibre orientations, and for $p=1.0001$}
\label{fig:H1_Error_p1}
\end{figure}

The following set of results show behaviour for various fibre orientations, and for values of $p=1.0001, 3$ and $10^4$.

Figure \ref{fig:H1_Error_p1} shows the $\mathcal{H}^1$-error convergence plots for all the formulations considered for $p=1.0001$.
Here $Q_1\_CG\_UI_\lambda$ and $Q_1\_CG\_UI_{\beta\lambda}$ show optimal convergence for any fibre direction at the superlinear rate $1.83$.
$Q_1\_CG$ and $Q_1\_CG\_UI_\beta$ show poor convergence.
\begin{figure}
\centering
\subfloat[$\hat{a} = \pi/4$]{\includegraphics[width=.5\columnwidth]{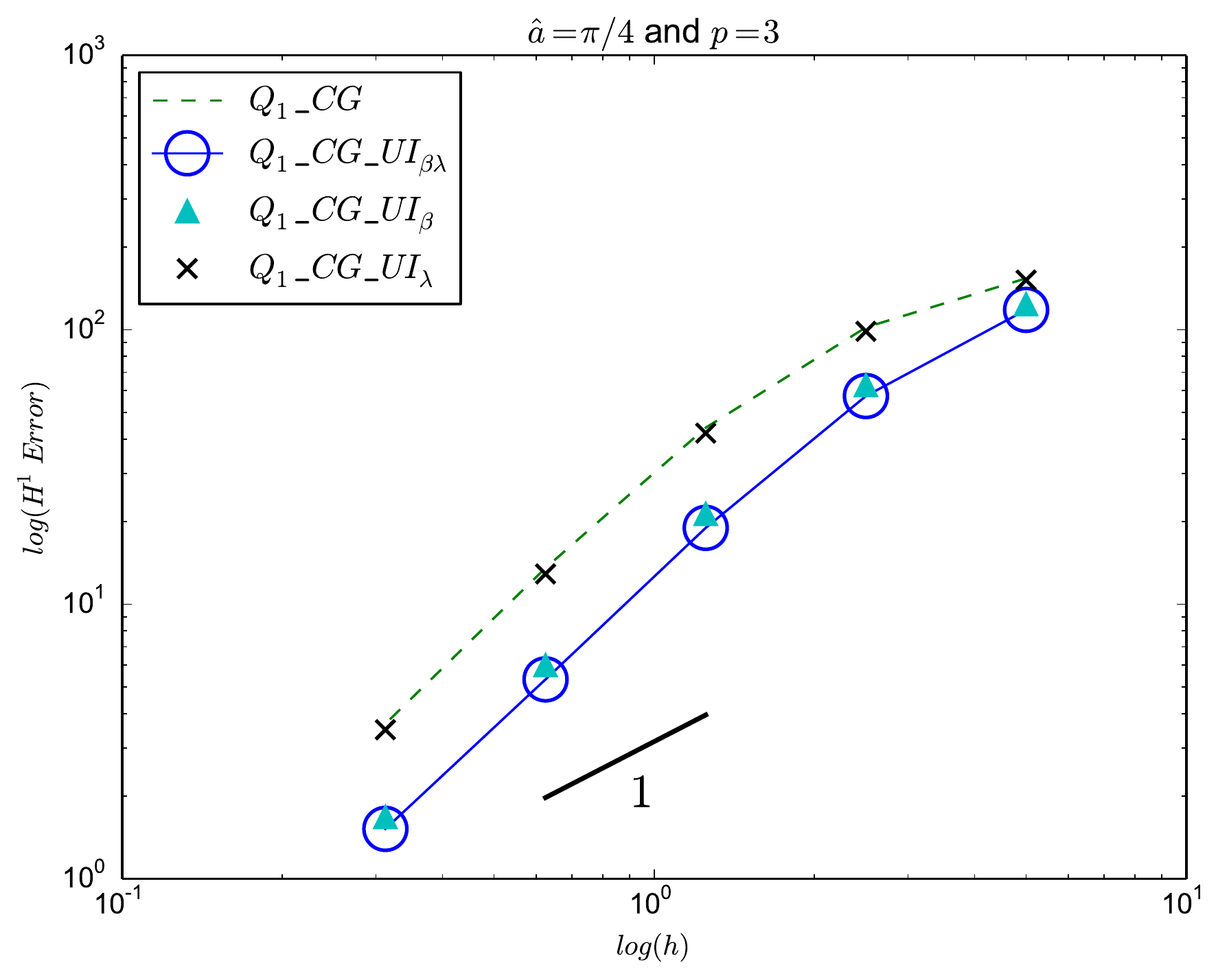}\label{fig:H1_p3_pi4}}
\subfloat[$\hat{a} = \pi/2$]{\includegraphics[width=.5\columnwidth]{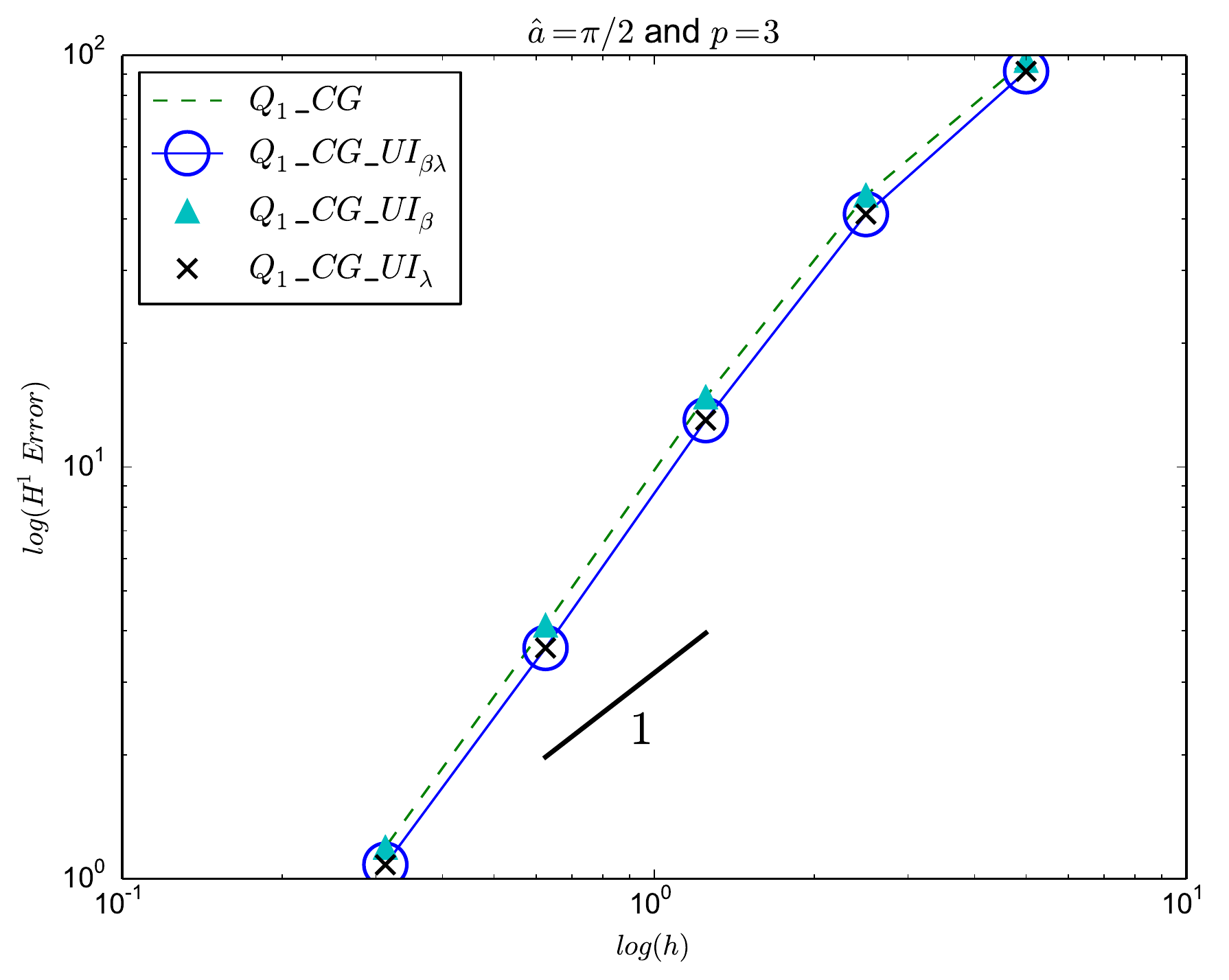}\label{fig:H1_p3_pi2}}\\
\subfloat[$\hat{a} = 3\pi/4$]{\includegraphics[width=.5\columnwidth]{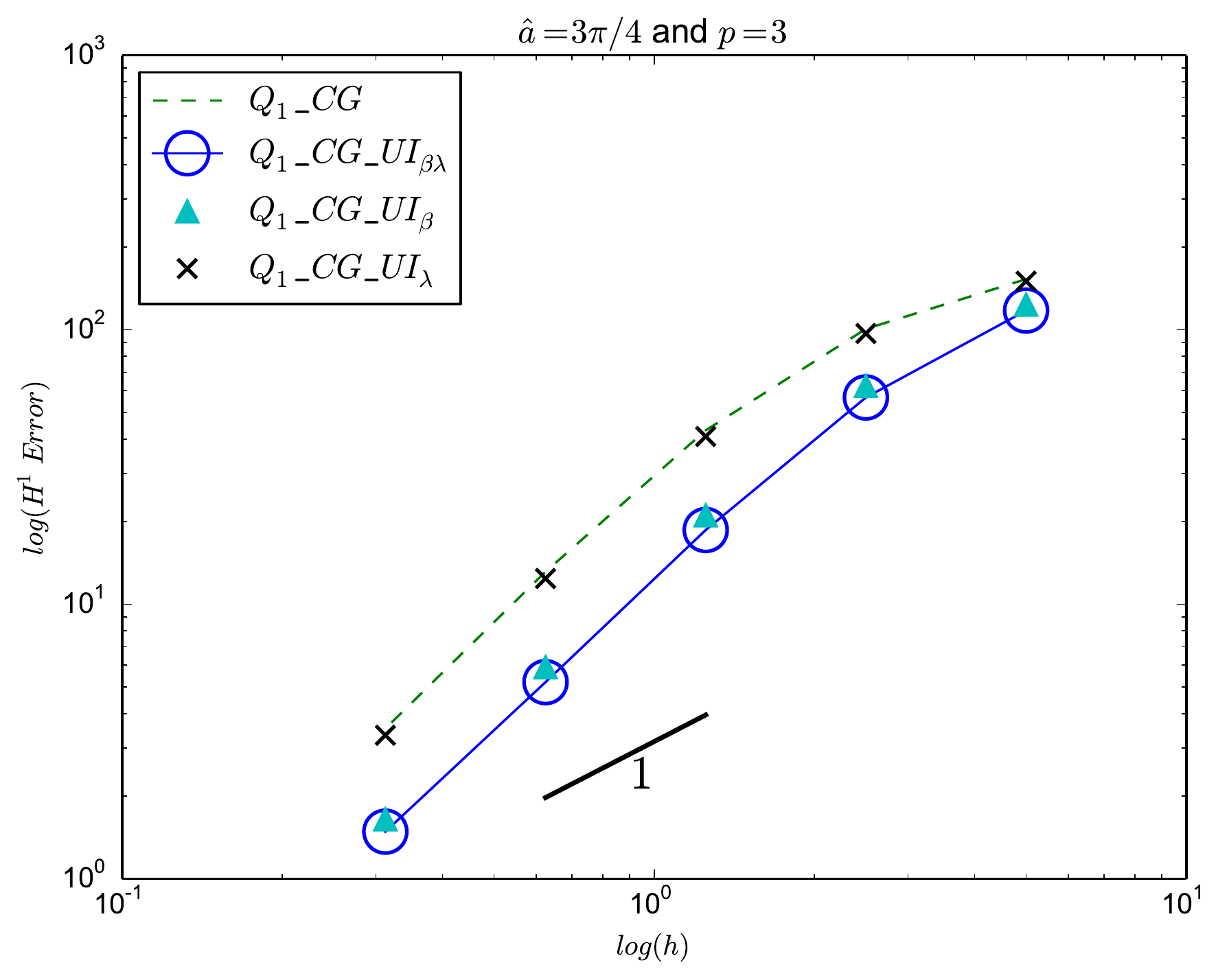}\label{fig:H1_p3_3pi4}}
\caption{Comparison of $\mathcal{H}^1$ errors for conforming and under-integrated elements on quadrilaterals, for different fibre orientations, and for $p=3$}
\label{fig:H1_Error_p3}
\end{figure}
Figure \ref{fig:H1_Error_p3} shows the $\mathcal{H}^1$-error convergence plots for all the formulations considered for $p=3$.
Here, all formulations at any fibre direction are superlinearly convergent at rate $1.7$.
\begin{figure}
\centering
\subfloat[$\hat{a} = \pi/4$]{\includegraphics[width=.5\columnwidth]{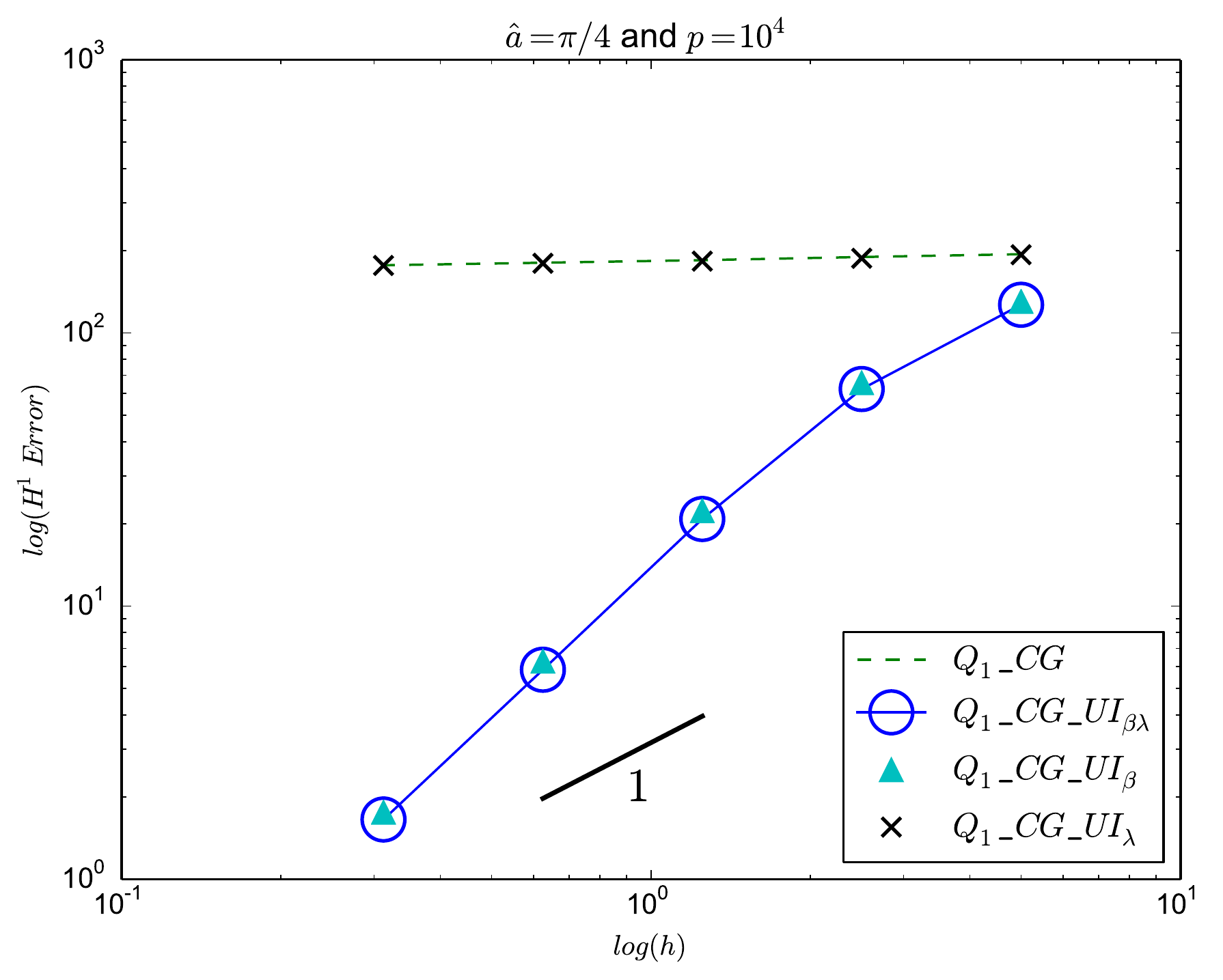}\label{fig:H1_p10e4_pi4}}
\subfloat[$\hat{a} = \pi/2$]{\includegraphics[width=.5\columnwidth]{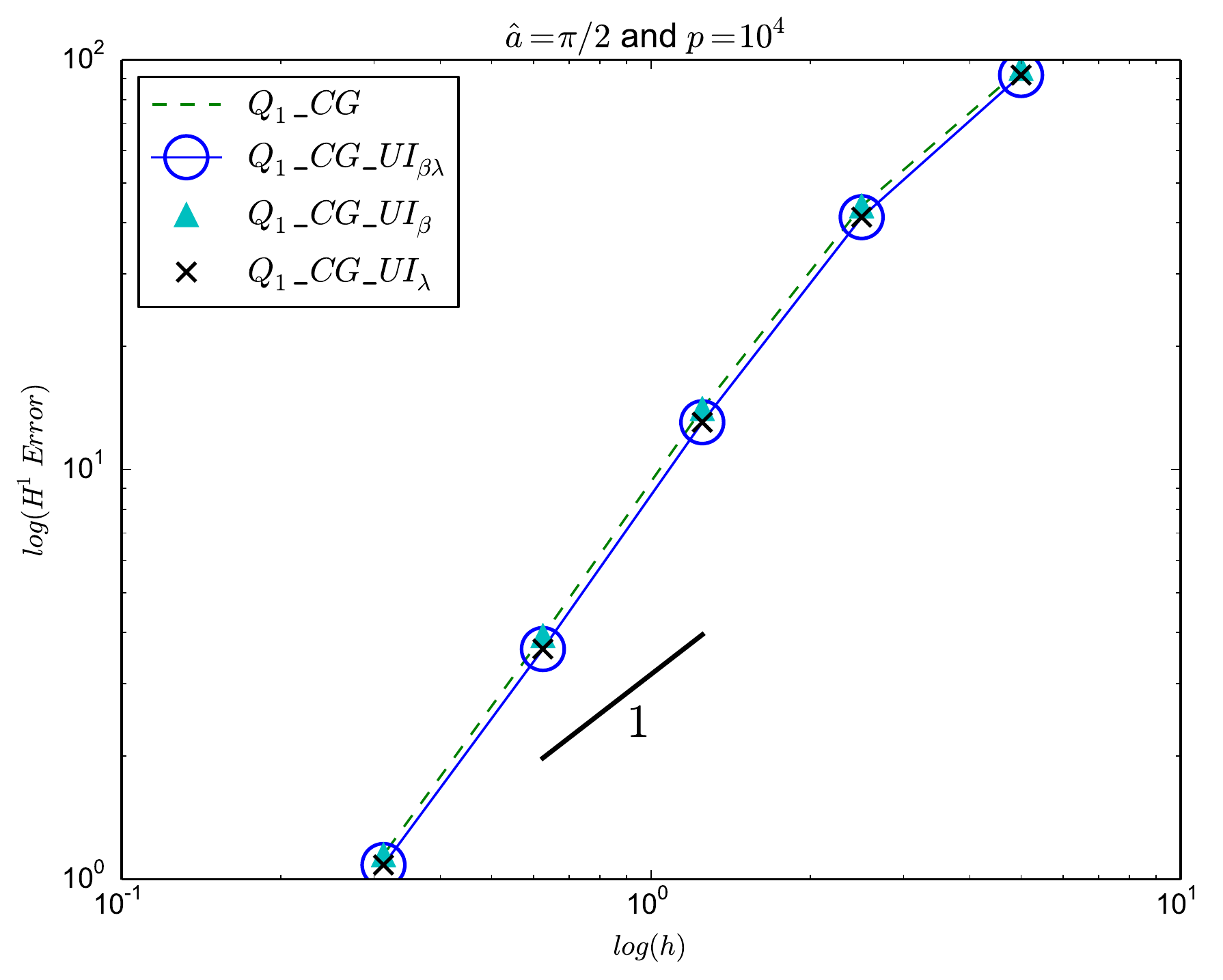}\label{fig:H1_p10e4_pi2}}\\
\subfloat[$\hat{a} = 3\pi/4$]{\includegraphics[width=.5\columnwidth]{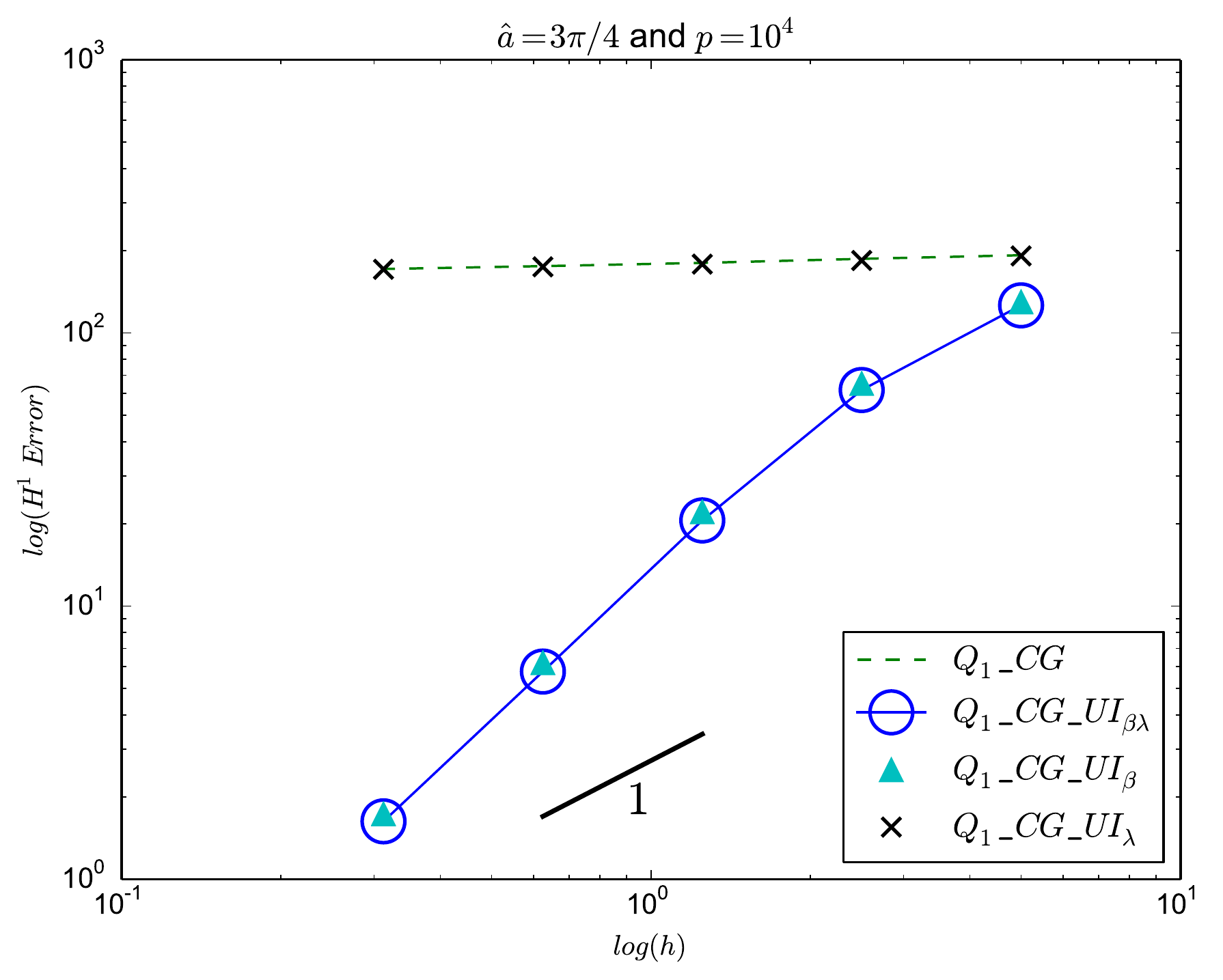}\label{fig:H1_p10e4_3pi4}}
\caption{Comparison of $\mathcal{H}^1$ errors for conforming and under-integrated elements on quadrilaterals, for different fibre orientations, and for $p=10^4$}
\label{fig:H1_Error_p10e4}
\end{figure}
Figure \ref{fig:H1_Error_p10e4} shows the $\mathcal{H}^1$-error convergence plots for all the formulations considered, for $p=10^4$.
Figure \ref{fig:H1_p10e4_pi4} and Figure \ref{fig:H1_p10e4_3pi4} show optimal convergence at the rate $1.86$ for $Q_1\_CG\_UI_\beta$ and $Q_1\_CG\_UI_{\beta\lambda}$,
and poor convergence for $Q_1\_CG$ and $Q_1\_CG\_UI_\lambda$.
In Figure \ref{fig:H1_p10e4_pi2} where the fibre direction is at an angle $\pi/2$, the error plots show convergence at the superlinear rate $1.64$ for $Q_1\_CG$.

%----------------------------------------------------------------------------------------
%		CONCLUSIONS
%----------------------------------------------------------------------------------------
\section{Conclusions}
The results in this work shed useful light on the relationship between well-posedness and anisotropy in the context of transverse isotropy, and show that, as far as computational approximations are concerned, transverse isotropy can have the effect of eliminating volumetric locking behaviour for moderate values of the anisotropy parameter, but away from its value of unity which would correspond to isotropy.
The use of under-integration in the extensional term follows the well-established approach used to circumvent volumetric locking; an extensive set of numerical results, over a range of measures of anisotropy and a range of fibre directions, has shown in this work that under-integration is robustly locking-free. 

The study and development of finite element formulations for incompressible and near-incompressible behaviour is extensive, with rigorous analyses supporting a range of stable and convergent element choices, whether with the use of mixed and enhanced formulations, equivalent under-integration schemes, or discontinuous Galerkin approaches, to mention some examples.
A corresponding framework is lacking in the case of near-inextensibility. The work reported here is intended to contribute to that framework.
Further investigations will involve the analysis and implementation of mixed and discontinuous Galerkin formulations. 

There have been a substantial number of investigations of near-inextensible behaviour under large-strain conditions.
It is clear from these studies, which were referred to in the Introduction, that the large-deformation problem presents features and challenges, some of which are absent in the small-strain regime.
It is intended to build on these investigations and the work reported here, in exploring theoretically and computationally the large-deformation problem.
%----------------------------------------------------------------------------------------
%		APPENDIX
%----------------------------------------------------------------------------------------
\begin{appendix}
\section*{Appendix}
For plane strain, the strain-stress relationship for a transversely isotropic material, with fibre direction $\bm{a} = \left(\begin{smallmatrix}a_1\\a_2\end{smallmatrix}\right)$, is
\[
\begin{pmatrix}
\varepsilon_{11}\\\varepsilon_{22}\\2\varepsilon_{12}
\end{pmatrix}
=
\begin{pmatrix}
\mathbb{S}_{11} & \mathbb{S}_{12} & \mathbb{S}_{13}\\
\mathbb{S}_{12} & \mathbb{S}_{22} & \mathbb{S}_{23}\\
\mathbb{S}_{13} & \mathbb{S}_{23} & \mathbb{S}_{33}
\end{pmatrix}
\begin{pmatrix}
\sigma_{11} \\ \sigma_{22} \\ \sigma_{12}
\end{pmatrix},
\]
where
\begin{small}
\begin{align*}
\mathbb{S}_{11} &= \dfrac{1}{det\mathbb{C}} \Big[ \big(\lambda+2\mu_t+2(\gamma+\alpha)a_2^2+\beta a_2^4\big)\big(\mu_t + \dfrac{\gamma}{2}+\beta a_1^2a_2^2\big) - \big((\alpha+\gamma)a_1a_2+\beta a_1a_2^3\big)^2\Big]\\
\mathbb{S}_{12} &= \dfrac{1}{det\mathbb{C}} \Big[ \big((\alpha+\gamma)a_1a_2+\beta a_1^3a_2\big)\big((\alpha+\gamma)a_1a_2+\beta a_1a_2^3\big) - \big(\lambda+\alpha+\beta a_1^2a_2^2\big)\big(\mu_t + \dfrac{\gamma}{2}+\beta a_1^2a_2^2\big)\Big]\\
\mathbb{S}_{13} &= \dfrac{1}{det\mathbb{C}} \Big[ \big(\lambda+\alpha+\beta a_1^2a_2^2\big)\big((\alpha+\gamma)a_1a_2+\beta a_1a_2^3\big) - \big((\alpha+\gamma)a_1a_2+\beta a_1^3a_2\big)\big(\lambda+2\mu_t+2(\gamma+\alpha)a_2^2+\beta a_2^4\big)\Big]\\
\mathbb{S}_{22} &= \dfrac{1}{det\mathbb{C}} \Big[ \big(\lambda+2\mu_t+2(\gamma+\alpha)a_1^2+\beta a_1^4\big)\big(\mu_t + \dfrac{\gamma}{2}+\beta a_1^2a_2^2 \big) - \big((\alpha+\gamma)a_1a_2+\beta a_1^3a_2\big)^2\Big]\\
\mathbb{S}_{23} &= \dfrac{1}{det\mathbb{C}} \Big[ \big(\lambda+\alpha+\beta a_1^2a_2^2\big)\big((\alpha+\gamma)a_1a_2+\beta a_1^3a_2\big) - \big(\lambda+2\mu_t+2(\gamma+\alpha)a_1^2+\beta a_1^4\big)\big((\alpha+\gamma)a_1a_2+\beta a_1a_2^3\big)\Big]\\
\mathbb{S}_{33} &= \dfrac{1}{det\mathbb{C}} \Big[ \big(\lambda+2\mu_t+2(\gamma+\alpha)a_1^2+\beta a_1^4\big)\big(\lambda+2\mu_t+2(\gamma+\alpha)a_2^2+\beta a_2^4\big) - \big(\lambda+\alpha+\beta a_1^2a_2^2\big)^2\Big]
\end{align*}
\end{small}
with
\begin{small}
\begin{align*}
\dfrac{1}{det\mathbb{C}} =& \big(\lambda+2\mu_t+2(\gamma+\alpha)a_1^2+\beta a_1^4\big) \Big[\big(\lambda+2\mu_t+2(\gamma+\alpha)a_2^2+\beta a_2^4\big) \big(\mu_t + \dfrac{\gamma}{2}+\beta a_1^2a_2^2\big)\\
			&\hspace{5.5cm} - \big((\alpha+\gamma)a_1a_2+\beta a_1a_2^3\big)^2\Big]\\
						 &- \big(\lambda+\alpha+\beta a_1^2a_2^2\big) \Big[\big(\lambda+\alpha+\beta a_1^2a_2^2\big) \big(\mu_t + \dfrac{\gamma}{2}+\beta a_1^2a_2^2\big)\\
			&\hspace{3.5cm} - \big((\alpha+\gamma)a_1a_2+\beta a_1^3a_2\big) \big((\alpha+\gamma)a_1a_2+\beta a_1a_2^3\big)\Big]\\
						 &+ \big((\alpha+\gamma)a_1a_2+\beta a_1^3a_2\big) \Big[\big(\lambda+\alpha+\beta a_1^2a_2^2\big) \big((\alpha+\gamma)a_1a_2+\beta a_1a_2^3\big)\\
			&\hspace{4.5cm} - \big((\alpha+\gamma)a_1a_2+\beta a_1^3a_2\big) \big(\lambda+2\mu_t+2(\gamma+\alpha)a_2^2+\beta a_2^4\big)\Big].
\end{align*}
\end{small}
\end{appendix}
%----------------------------------------------------------------------------------------
% ACKNOWLEDGEMENT
%----------------------------------------------------------------------------------------
\section*{Acknowledgements}
The authors acknowledge with thanks the support for this work by the National Research Foundation, through the South African Research Chair in Computational Mechanics.
%----------------------------------------------------------------------------------------
%	BIBLIOGRAPHY
%----------------------------------------------------------------------------------------
%\renewcommand{\refname}{\spacedlowsmallcaps{References}} % For modifying the bibliography heading
\bibliographystyle{unsrt}
\bibliography{TI.bib} % The file containing the bibliography
%----------------------------------------------------------------------------------------

\end{document}